\documentclass{article}

\usepackage{amsmath,amsfonts,latexsym,amssymb,amsthm,url}

\usepackage[toc]{appendix}

\usepackage[utf8]{inputenc}

\makeatletter
\let\@fnsymbol\@alph
\makeatother

\newcommand{\C}{{\mathbb C}}
\newcommand{\F}{{\mathbb F}}
\renewcommand{\H}{{\mathbb H}}
\newcommand{\Q}{{\mathbb Q}}

\newcommand{\Z}{{\mathbb Z}}

\DeclareMathOperator{\lcm}{lcm}
\DeclareMathOperator{\Gal}{Gal}
\DeclareMathOperator{\height}{h}
\DeclareMathOperator{\id}{id}
\renewcommand{\Im}{\mathrm{Im}}
\newcommand{\Mu}{\mathrm{M}}
\renewcommand{\Re}{\mathrm{Re}}
\newcommand{\SL}{\mathrm{SL}}
\newcommand{\tho}{{\mathrm{th}}}

\newcommand{\calC}{\mathcal{C}}

\newcommand{\calF}{\mathcal{F}}

\newcommand{\calO}{\mathcal{O}}

\newcommand{\calQ}{\mathcal{Q}}

\newcommand{\bfm}{\mathbf{m}}

\newcommand\eps\varepsilon
\newcommand\ph\varphi

\newcommand{\tilgamma}{\tilde{\gamma}}
\newcommand{\tileta}{\tilde{\eta}}

\newcommand{\gerp}{\mathfrak{p}}

\newtheorem{theorem}{Theorem}[section]

\newtheorem{proposition}[theorem]{Proposition}

\newtheorem{corollary}[theorem]{Corollary}

\newtheorem{lemma}[theorem]{Lemma}

\newtheorem{remark}[theorem]{Remark}

\numberwithin{equation}{section}

\title{Effective multiplicative independence of~$3$ singular moduli}

\author{Yuri Bilu\footnote{Supported by SPARC Project P445 (India).}\textsuperscript{\, ,}\footnote{Supported by ANR project JINVARIANT.}, \ Sanoli Gun\textsuperscript{a},\ and Emanuele Tron\textsuperscript{b}}

\date{}

\makeatletter

\renewcommand*\l@section[2]{
	\ifnum \c@tocdepth >\z@
	\addpenalty\@secpenalty
	\addvspace{0.2em \@plus\p@}
	\setlength\@tempdima{1.5em}
	\begingroup
	\parindent \z@ \rightskip \@pnumwidth
	\parfillskip -\@pnumwidth
	\leavevmode \bfseries
	\advance\leftskip\@tempdima
	\hskip -\leftskip
	#1\nobreak\hfil \nobreak\hb@xt@\@pnumwidth{\hss #2}\par
	\endgroup
	\fi}

\makeatother

\begin{document}
	
	\hfuzz 5pt

	\maketitle

\begin{abstract}
Pila and Tsimerman  proved in 2017 that for every~$k$ there exists at most finitely many $k$-tuples $(x_1,\ldots, x_k)$ of non-zero singular moduli such
that $x_1, \ldots,x_k$ are multiplicatively dependent, but any proper subset of them is multiplicatively independent.
The proof was non-effective, using Siegel's lower bound for the class 	numbers. In 2019 Riffaut obtained an effective version of this result for ${k=2}$. Moreover, he determined all the instances of ${x^my^n\in \Q^\times}$, where $x,y$ are distinct singular moduli and $m,n$ are non-zero integers. In this article we obtain a similar result for ${k=3}$. We show that ${x^my^nz^r\in \Q^\times}$ (where $x,y,z$ are distinct singular moduli and $m,n,r$ non-zero integers) implies that the discriminants of $x,y,z$ do not exceed $10^{10}$. 
\end{abstract}

	{\footnotesize
		
		\tableofcontents
		
	}

	\section{Introduction}
	
	A \textit{singular modulus} is the $j$-invariant of an elliptic curve with complex multiplication. Given a singular modulus~$x$ we denote by $\Delta_x$ the discriminant of the associated imaginary quadratic order. 
	We denote by $h(\Delta)$ the class number of the imaginary quadratic order of discriminant~$\Delta$.
	
	Recall that two singular moduli~$x$ and~$y$ are conjugate over~$\Q$ if and only if ${\Delta_x=\Delta_y}$, and that there are $h(\Delta)$ singular moduli of a given discriminant~$\Delta$.   In particular, ${[\Q(x):\Q]=h(\Delta_x)}$. For all details see, for instance, \cite[\S~7 and \S~11]{Co13}.
	
%\subsection{History}	
	
	There has been much work on diophantine properties of singular moduli in recent years. In particular, studying algebraic equations where the unknowns are singular moduli \cite{ABP15,BLM17,BLP16} is interesting by virtue of its connection with the André--Oort property for affine space \cite{BMZ13,Ku12,Ku13}.
	
	Pila and Tsimerman~\cite{PT17} proved that for every~$k$ there exists at most finitely many $k$-tuples $(x_1,\ldots, x_k)$ of  non-zero singular moduli such that
	$x_1, \ldots,x_k$ are multiplicatively dependent, but any proper subset of them is multiplicatively independent. Their argument is fundamentally non-effective.  
	
	Riffaut  \cite[Th.~1.7]{Ri19} gave an effective (in fact, totally explicit) version of the theorem of Pila and Tsimerman in the case ${k=2}$. In fact, he did more: he classified all cases when  ${x^my^n\in \Q^\times}$, where $x,y$ are singular moduli and $m,n$ non-zero integers. 

%\subsection{Results}	
	
	In the present article we obtain an effective result for ${k=3}$. Moreover, like Riffaut did, we prove a stronger statement: we bound explicitly discriminants of singular moduli $x,y,z$ such that ${x^my^nz^r\in \Q^\times}$ for some non-zero integers $m,n,r$. Our bound is as follows.
	
	\begin{theorem}
		\label{thpizthree}
		Let $x,y,z$ be distinct non-zero singular moduli and ${m,n,r}$ non-zero integers. Assume that ${x^my^nz^r\in \Q^\times}$. Then 
		$$
		\max\{|\Delta_x|,|\Delta_y|,|\Delta_z|\}< 10^{10}.
		$$

	\end{theorem}
	
	The special case ${m=n=r}$ has been recently settled by Fowler~\cite{Fo20,Fo23}. 
	
Note that there do exist triples of distinct singular moduli $x,y,z$ such that ${x^my^nz^r\in \Q^\times}$ for some non-zero ${m,n,r\in \Z}$. There are three types of  currently known examples.

\begin{description}
\item[Rational type]
Take distinct $x,y,z$ such that 
$$
h(\Delta_x)=h(\Delta_y)=h(\Delta_z)=1,\qquad \Delta_x,\Delta_y,\Delta_z\ne -3. 
$$ 
In this case ${x,y,z\in \Q^\times}$ and ${x^my^nz^r\in \Q^\times}$ for any choice of $m,n,r$. Pila and Tsimerman \cite[Example 6.2]{PT17} even found an example of ${x^my^nz^r=1}$:
$$
(2^63^3)^{10}(-2^{15})^6(-2^{15}3^3)^{-10}=1, 
$$
the corresponding discriminants being $-4$, $-11$ and $-19$. 

\item[Quadratic type]
Take distinct $x,y,z$ such that 
$$
h(\Delta_x)=1,\qquad \Delta_x\ne -3, \qquad \Delta_y=\Delta_z, \qquad h(\Delta_y)=h(\Delta_z)=2. 
$$ 
In this case ${x\in \Q^\times}$ and ${y,z}$ are of degree~$2$, conjugate over~$\Q$. Hence ${x^my^nz^n\in \Q^\times}$ for any choice of $m,n$.

\item[Cubic type]
Take distinct $x,y,z$ such that 
$$
\Delta_x=\Delta_y=\Delta_z,\qquad  h(\Delta_x)=h(\Delta_y)=h(\Delta_z)=3. 
$$ 
In this case  ${x,y,z}$ are of degree~$3$, forming a full Galois orbit over~$\Q$.  Hence ${xyz\in \Q^\times}$. 
\end{description}

We believe that, up to permuting $x,y,z$, there are no other examples, but to justify it, one needs to improve on the numerical bound $10^{10}$ in Theorem~\ref{thpizthree}. 	
	
	The proof of Theorem~\ref{thpizthree} relies on the following result, which is a
	partial common generalization (for big discriminants) of \cite[Theorem~1.7]{Ri19} and \cite[Theorem~1.3]{FR18}. 
	
	\begin{theorem}
		\label{thprimel}
		Let $x,y$ be distinct non-zero singular moduli and $m,n$ non-zero integers. Assume that 
		\begin{equation}\label{condthm}
		\max\{|\Delta_x|,|\Delta_y|\}\ge 10^8. 
		\end{equation}
		Then ${[\Q(x,y):\Q(x^my^n)]\le 2}$. More precisely, we have either 
		\begin{equation}
		\label{esamefield}
		\Q(x^my^n)=\Q(x,y)
		\end{equation}
		or 
		\begin{equation}
		\label{esubfield}
		\Delta_x=\Delta_y, \qquad m=n, \qquad [\Q(x,y):\Q(x^my^m)]= 2.
		\end{equation}
		Moreover, in the latter case~$x$ and~$y$ are conjugate over the field ${\Q(x^my^m)}$. 
		
		If $\{ \Delta_x,\Delta_y \}$ is not of the form $\{\Delta,4\Delta\}$ for some $\Delta \equiv 1 \bmod 8$, then condition \eqref{condthm} can be relaxed to \begin{equation}
		\max\{|\Delta_x|,|\Delta_y|\}\ge  10^6. 
		\end{equation}
	\end{theorem}

%\subsection{Examples}

\paragraph{Plan of the article}
	
	In Section~\ref{sprelimi}  we collect basic fact about singular moduli to be used throughout the article. In Section~\ref{slinear} we establish our principal tool: a linear relation between the exponents ${m_1, \ldots, m_k}$ stemming from the multiplicative relation ${x_1^{m_1}\cdots x_k^{m_k}=1}$. Theorems~\ref{thprimel} and~\ref{thpizthree} are proved in Sections~\ref{sprimel} and~\ref{spizthree}, respectively. 
	
	\paragraph{Acknowledgments} We thank Guy Fowler and the anonymous referee for many  comments that helped us to correct inaccuracies and improve the presentation.  We also thank Francesco Amoroso, Margaret Bilu, Igor Rapinchuk and Anatoly Vorobey for useful suggestions. Finally, we thank Bill Allombert and Amalia Pizarro, who allowed us to borrow Proposition~\ref{pbilly} from~\cite{ABP14}.  
	
	All calculations were performed using \textsf{PARI}~\cite{pari}.  We thank Bill Allombert and Karim Belabas for the \textsf{PARI} tutorial. The reader may consult 
	\begin{center}
		\url{https://github.com/yuribilu/multiplicative} 
	\end{center}
	to view the \textsf{PARI} scripts used for this article.
	
	\subsection{Notation and conventions}
	\label{ssnota}
	We denote by~$\H$ the Poincaré half-plane, and by~$\calF$ the standard fundamental domain for the action of the modular group; that is, the open hyperbolic triangle with vertices 
	$$
	\zeta_6=\frac{1+\sqrt{-3}}{2}, \quad \zeta_3=\frac{-1+\sqrt{-3}}{2}, \quad i\infty,
	$$
	together with the hyperbolic geodesics $[i,\zeta_6]$ and ${[\zeta_6,i\infty]}$.
	
	We denote by $\log$ the principal branch of the complex logarithm:
	$$
	-\pi < \arg\log z\le \pi \qquad (z\in \C^\times). 
	$$
	
	We use ${O_1(\cdot)}$ as a quantitative version of the ${O(\cdot)}$ notation: ${A=O_1(B)}$ means that ${|A|\le B}$. 
	
	We write the Galois action exponentially: ${x\mapsto x^\sigma}$. In particular, it is a right action: ${x^{(\sigma_1\sigma_2)}=(x^{\sigma_1})^{\sigma_2}}$. Most of the Galois groups occuring in this article are abelian, so this is not relevant, but in the few cases where the group is not abelian one must be vigilant. 
	
	Let~$R$ be a commutative ring, and ${a\in R}$. When this does not lead to confusion, we write $R/a$ instead of $R/aR$.

We denote by $\calC_m$ the cyclic group of order~$m$. 	 
	
	One point on references: item~Y of Proposition~X is quoted as Proposition~X:Y.

	\section{Class numbers, denominators, isogenies}
	\label{sprelimi}

	Unless the contrary is stated explicitly,  the letter~$\Delta$ stands for an \textit{imaginary quadratic discriminant}, that is, ${\Delta<0}$ and  ${\Delta\equiv 0,1\bmod 4}$. 
	
	We denote by $\calO_\Delta$ the imaginary quadratic order of discriminant~$\Delta$, that is, ${\calO_\Delta =\Z[(\Delta+\sqrt\Delta)/2]}$. Then ${\Delta=Df^2}$, where~$D$ is the  discriminant of the number field ${K=\Q(\sqrt\Delta)}$ (called the \textit{fundamental discriminant} of~$\Delta$) and ${f=[\calO_D:\calO_\Delta]}$ is called the \textit{conductor} of~$\Delta$.
	
	We denote by $h(\Delta)$ the class  number of~$\calO_\Delta$.
	
	Given a singular modulus~$x$, we denote $\Delta_x$ the discriminant of the associated CM order, and  we write ${\Delta_x=D_xf_x^2}$ with~$D_x$ the fundamental discriminant and~$f_x$ the conductor. 
	Furthermore, we denote by~$K_x$ the associated imaginary quadratic field
	$$
	K_x=
	\Q(\sqrt{D_x})=\Q(\sqrt{\Delta_x}).
	$$
	We will call~$K_x$ the \textit{CM-field} of the singular modulus~$x$.

	It is known (see, for instance, §11 in~\cite{Co13}) that a singular modulus~$x$ is an algebraic integer of degree $h(\Delta_x)$, and that there are exactly $h(\Delta)$  singular moduli of given discriminant~$\Delta$, which form a full Galois orbit over~$\Q$. 
	
	\subsection{Class numbers and class groups}
	
	For a discriminant~$\Delta$ and a positive integer~$\ell$ set 
	\begin{equation}
	\label{epsild} \Psi(\ell,\Delta)=\ell\prod_{p\mid \ell}\left(1-\frac{(\Delta/p)}{p}\right),  
	\end{equation}
where $(\Delta/p)$ denotes the Kronecker symbol. 	It is useful to note that 
	\begin{equation}
	\label{elowereuler}
	\Psi(\ell,\Delta)\ge \ph(\ell),
	\end{equation}
	where ${\ph(\cdot)}$ is Euler's totient function. Note also the multiplicativity relation
	\begin{equation}
	\label{emult}
	\Psi(\ell_1\ell_2, \Delta)= \Psi(\ell_2, \Delta\ell_1^2) \Psi(\ell_1,\Delta). 
	\end{equation}

	Recall the \textit{class number formula}:  
	\begin{equation}
	\label{eclnfr}
	h(\Delta\ell^2) = \frac{1}{[\calO_\Delta^\times:\calO_{\Delta\ell^2}^\times]}h(\Delta)\Psi(\ell,\Delta) , 
	\end{equation}
	see \cite[Theorem~7.24]{Co13}. Note that in~\cite{Co13} it is proved only in the case when ${\Delta=D}$ is a fundamental discriminant. However, the general case easily follows from the case ${\Delta=D}$ using the multiplicativity relation~\eqref{emult}. Note also that 
\begin{equation}
\label{efactor}
	[\calO_\Delta^\times:\calO_{\Delta\ell^2}^\times]=
	\begin{cases}
	3,& \Delta=-3,\ell>1,\\
	2,& \Delta=-4,\ell>1,\\
	1,& \Delta\ne -3,-4. 
	\end{cases}
\end{equation}
	
\subsubsection{Discriminants with small class number}	
	
	Watkins~\cite{Wa04} classified fundamental discriminants~$D$ with ${h(D)\le 100}$. In particular, he proved that such discriminants do not exceed $2383747$ in absolute value. It turns out that the same upper bound holds true for all discriminants, not only for fundamental ones.

	\begin{proposition}
		\label{pwatki}
		Let~$\Delta$ be a negative discriminant with ${h(\Delta) \le 100}$. Then we have ${|\Delta|\le 2383747}$. If $h(\Delta) \leq 64$, then $|\Delta| \leq 991027$.
	\end{proposition}

\begin{remark}
As Guy Fowler informed us, Janis Klaise obtained the same result, with a similar proof, in his 2012 Master Thesis~\cite{KL12}. Apparently, his manuscript had never been published. 
\end{remark}	
	
	\begin{proof} 
		Given a positive integer~$n$, denote 
		$$
		D_{\max}(n):=\max\{|D|: D \text{ fundamental}\mathbin{,} h(D)\le n\}
		$$
		the biggest absolute value of a \textit{fundamental} discriminant~$D$ with ${h(D)\le n}$; the values of $D_{\max}$ for arguments up to $100$ can be found in Watkins \cite[Table~4]{Wa04}. For the reader's convenience, we give in Table~\ref{tadmax} the  $D_{\max}$ of the arguments occurring in equation~\eqref{edeltale} below. 
		
		\begin{table}
			\caption{Values of $D_{\max}(n)$ for~$n$ of the form $\lfloor 100/\varphi(f)\rfloor$}
			\label{tadmax}
			{\scriptsize
				$$
				\begin{array}{rr|ccccccccc}
				&n&1&2&3&4&5&6&7&8&10\\
				&\max f &420&210&120&90 & 66 & 60 & 42 & 42 & 30\\
				&D_{\max}(n)&163&427&907&1555&2683&3763&5923&6307&13843\\
				{ \ }\\
				&n&12&16&25&50&100\\
				&\max f&30 & 18 & 12 &6& 2  \\
				&D_{\max}(n)& 17803 & 34483& 111763 & 462883 & 2383747
				\end{array}
				$$
				
				\textbf{Explanation.\quad}
				The first row contains all positive integers~$n$ of the form $\lfloor 100/\varphi(f)\rfloor$ for some positive integer~$f$. In the second row, for each~$n$ we give the biggest~$f$ with the property ${ 100/\varphi(f)\ge n}$. Finally, in the third row we display $D_{\max}(n)$. 
			}
		\end{table}

		Now let ${\Delta=Df^2}$ be such that ${h(\Delta) \le 100}$. Using the class number formula~\eqref{eclnfr} (applied with~$D$ as~$\Delta$ and with~$f$ as~$\ell$), and the bound~\eqref{elowereuler} we get 
		\begin{equation}
		\label{eintermed} 
		h(D) \varphi(f) \leq 100 [\calO_D^\times:\calO_\Delta^\times]. 
		\end{equation}
		
		If $D=-3$ or $-4$ then this implies $\varphi(f) \leq 300$: the largest such $f$ is $f=1260$, so that in this case ${|\Delta|\leq 6350400}$.
		
		If $D \neq -3,-4$ then we find from \eqref{eintermed} that $h(D) \leq 100/\varphi(f)$, and hence 
		\begin{equation}
		\label{edeltale} 
		|\Delta|=f^2|D| \leq f^2  D_{\max}(\left\lfloor 100/\varphi(f)\right\rfloor).  
		\end{equation}
		Plugging in the values from Table~\ref{tadmax}, we find that the maximum of the right-hand side is attained for $f=420$ and is equal to $28753200$. This proves that ${|\Delta|\le 28753200}$. 
		
		To complete the proof, we run a \textsf{PARI} script computing the class numbers of all~$\Delta$ with ${|\Delta|\le 28753200}$. It confirms that the biggest~$\Delta$ with ${h(\Delta)\le 100}$ is ${-2383747}$, and the biggest~$\Delta$ with ${h(\Delta)\le 64}$ is ${-991027}$. %The total running time was 50 minutes on a personal computer. 
	\end{proof}
	
\subsubsection{The $2$-rank}
\label{ssstwor}
	
	Given a finite abelian group~$G$ and a prime number~$p$, we call the \textit{$p$-rank} of~$G$, denoted $\rho_p(G)$, the dimension of the $\F_p$-vector space ${G/G^p}$. If~$\Delta$ is a discriminant then we denote by $\rho_p(\Delta)$ the $p$-rank of its class group. 
	
	The $2$-rank of   a discriminant was determined by Gauss, see \cite[Proposition 3.11 and Theorem~3.15]{Co13}. As usual, we denote by $\omega(n)$  the number of distinct prime divisors of a non-zero integer~$n$.	  %We will express this result in the language of ring class fields. 
%	We do not need the full strength of the Gauss Theorem, but only the following consequence. 
\begin{proposition}
		\label{pgauss}
		%Let~$L$ be the ring class field of  discriminant~$\Delta$ and let ${K=\Q(\sqrt{\Delta})}$ be its CM-field. 
		Let~$\Delta$ be a discriminant. Then
$$
\rho_2(\Delta)=
\begin{cases}
\omega(\Delta) -1, & \Delta\equiv 1 \bmod 4,\\
\omega(\Delta)-2, & \Delta \equiv 4 \bmod 16, \\
\omega(\Delta) -1, & \Delta \equiv 8,12 \bmod 16,\\
\omega(\Delta) -1,& \Delta\equiv 16 \bmod 32, \\
\omega(\Delta),  & \Delta\equiv 0\bmod 32. 
\end{cases}
$$			
In particular,
${\rho_2(\Delta) \in \{\omega(\Delta), \omega(\Delta)-1, \omega(\Delta)-2\}}$. If~$D$ is a fundamental discriminant, then 
		${\rho_2(D)=\omega(D)-1}$. 		
\end{proposition}

\subsection{Ring class fields}

	If~$x$ is a singular modulus with discriminant ${\Delta=Df^2}$ and ${K=\Q(\sqrt D)}$ is its CM-field, then $K(x)$ is an abelian extension of~$K$ such that $\Gal(K(x)/K)$ is isomorphic to to the class group of~$\Delta$; in particular,   ${[K(x):K]=h(\Delta)}$, and the singular moduli of discriminant~$\Delta$ form a full Galois orbit over~$K$ as well.

	This leads to the useful notion of \textit{ring class field}. Given an imaginary quadratic field~$K$ of discriminant~$D$ and a positive integer~$f$, the \textit{ring class field} of~$K$ of conductor~$f$, denoted $K[f]$, is, by definition, $K(x)$, where~$x$ is some singular modulus of discriminant $Df^2$. It does not depend on the particular choice of~$x$ and is an abelian extension of~$K$. 
	
	Proofs of the statements above can be found, for instance, in §9--11 of~\cite{Co13}.

	\bigskip
	
	The following properties will be systematically used. 
	
	\begin{proposition}
		\label{pdihedral}
		Let~$K$ be an imaginary quadratic field and~$L$ a ring class field of~$K$. Denote
		\begin{equation}
		\label{egh}
		G=\Gal(L/\Q), \qquad H=\Gal(L/K).
		\end{equation}
		(As we have just seen,~$H$ is an abelian group.) 
		Then we have the following. 
		\begin{itemize}
			
			\item
			every element of ${G\smallsetminus H}$ is of order~$2$;

			\item
			if ${\gamma\in G\smallsetminus H}$ and ${\eta\in H}$ then ${\gamma\eta\gamma=\eta^{-1}}$. 
		\end{itemize}
		
		In particular, any element of ${G\smallsetminus H}$ does not commute with any element of~$G$ of order bigger than~$2$. 
	\end{proposition}
	
	For the proof see, for instance, \cite[Lemma~9.3]{Co13}. 
	
	\begin{proposition}
		\label{pcompo}
		Let~$K$ be an imaginary quadratic field of discriminant~$D$, and $\ell,m$ positive integers. 
		\begin{enumerate}
			\item
			Assume that either ${D\ne -3,-4}$ or ${\gcd(\ell,m)>1}$. Then the compositum $K[\ell]K[m]$ is equal to $K[\lcm(\ell,m)]$.
			
			\item
			Assume that ${D=-3}$ and ${\gcd(\ell, m)=1}$.   Then $K[\ell]K[m]$ is either equal to $K[\lcm(\ell,m)]$ or is a subfield of $K[\lcm(\ell,m)]$ of degree~$3$.  
			
			\item
			Assume that ${D=-4}$ and ${\gcd(\ell, m)=1}$.  Then $K[\ell]K[m]$ is either equal to $K[\lcm(\ell,m)]$ or is a subfield of $K[\lcm(\ell,m)]$ of degree~$2$.  
		\end{enumerate}

	\end{proposition}
	
	For the proof see, for instance, \cite[Proposition~3.1]{ABP15}

	\subsubsection{Two-elementary subfields of ring class fields}
	
	We call a group \textit{$2$-elementary} if all its elements are of order dividing~$2$. A finite $2$-elementary group is a product of cyclic groups of order~$2$. Let~$K$ be a field and~$L$ a finite  extension of~$K$; we say that~$L$ is $2$-elementary over~$K$ if~$L$ is Galois over~$K$, with $2$-elementary Galois group.   We call a number field $2$-elementary if it is  $2$-elementary over~$\Q$.

	The following is well-known, but we include the proof for the reader's convenience. 
	
	\begin{proposition}
		\label{pabelian}
		Let~$F$ be a number field abelian over~$\Q$ and contained in some ring class field. 
		Then~$F$ is $2$-elementary. 
	\end{proposition}
	
	\begin{proof}
		This is an easy consequence of Proposition~\ref{pdihedral}. 
		Let~$K$ be an imaginary quadratic field such that its ring class field, denoted~$L$, contains~$F$. 
		We use the notation of~\eqref{egh}. 
	
		For ${\gamma \in G}$ let ${\tilgamma \in \Gal(F/\Q)}$ denote the restriction of~$\gamma$ to~$F$. Each element of $\Gal(F/\Q)$ is a restriction of either  some ${\gamma\in G\smallsetminus H}$ or some ${\eta\in H}$. In the former case ${\tilgamma^2=1}$ because ${\gamma^2=1}$. Now consider~$\tileta$ for some ${\eta\in H}$. Pick ${\gamma \in G\smallsetminus H}$. Then ${\tilgamma\tileta\tilgamma =\tileta^{-1}}$. But ${\Gal(F/\Q)}$ is abelian, which implies that ${\tilgamma\tileta\tilgamma =\tilgamma^2\tileta =\tileta}$. Hence ${\tileta^2=1}$ as well. Thus, every element of ${\Gal(F/\Q)}$ is of order dividing~$2$, as wanted.
	\end{proof}
	
	The only positive integers~$m$ such that the multiplicative group ${(\Z/m\Z)^\times}$ is $2$-elementary are the divisors of~$24$. Hence we have the following corollary.
	
	\begin{corollary}
		\label{crofoneinrcf}
		The group of roots of unity in a ring class field is of order dividing~$24$. 
	\end{corollary}

	Another important case of $2$-elementary fields is the intersection ${\Q(x)\cap\Q(y)}$, where~$x$ and~$y$ singular moduli with distinct fundamental discriminants. This is known since ages (see, for instances, the articles of André~\cite{An98} or Edixhoven~\cite{Ed98}), but we again include a proof for the reader's convenience.

	\begin{proposition}
		\label{pintersect}
		Let~$x$ and~$y$ be singular moduli with distinct fundamental discriminants: ${D_x\ne D_y}$. Then the field ${\Q(x)\cap \Q(y)}$ is $2$-elementary. In particular, if ${\Q(x)\subset \Q(y)}$ then $\Q(x)$ is $2$-elementary.
	\end{proposition}

	\begin{proof}
		It suffices to prove that the field ${\Q(x)\cap \Q(y)}$ is  abelian: Proposition~\ref{pabelian} will then complete the job.
		
		Recall that we let ${K_x=\Q(\sqrt{D_x})}$ denote the CM field for~$x$. We will denote $K_{xy}$ the compositum of~$K_x$ and~$K_y$, that is, the field $\Q(\sqrt{D_x},\sqrt{D_y})$. Furthermore, we define %let~$M$ denote the compositum and~$L$  the intersection of the fields  $K_{xy}(x)$ and $K_{xy}(y)$:
		$$
		M=K_{xy}(x,y), \qquad L= K_{xy}(x)\cap K_{xy}(y). 
		$$
		It suffices to prove that~$L$ is abelian, because ${L\supset \Q(x)\cap \Q(y)}$. To start with, let us prove that~$L$ is  $2$-elementary over the field $K_{xy}$.

		Since ${K_x\ne K_y}$, there exists  ${\iota \in \Gal(M/\Q)}$  such that 
		$$
		\iota\vert_{K_x}=\id \quad\text{and}\quad  \iota\vert_{K_y}\ne \id. 
		$$
		Proposition~\ref{pdihedral} implies that  for  ${\eta \in \Gal(M/K_{xy})}$ we have 
		$$
		\iota^{-1}\eta \iota\vert_{K_x(x)}= \eta\vert_{K_x(x)} \quad\text{and}\quad \iota^{-1}\eta \iota\vert_{K_y(y)}= \eta^{-1}\vert_{K_y(y)}. 
		$$
		We also have ${\eta\vert_{K_{xy}}=\id}$ by the choice of~$\eta$. Hence ${\eta\vert_L=\eta^{-1}\vert_L}$. Since every element of $\Gal(L/K_{xy})$ is a restriction to~$L$ of some ${\eta \in \Gal(M/K_{xy})}$, this proves that the Galois group $\Gal(L/K_{xy})$ is $2$-elementary, as wanted.

		To complete the proof, we must show that~$L$ is abelian over~$\Q$. Clearly,~$L$ is Galois over~$\Q$, being the intersection of two Galois extension. We have to show that ${\Gal(K_{xy}/\Q)}$ acts trivially on ${\Gal(L/K_{xy})}$. This means proving the following: for every ${\eta,\gamma \in \Gal(M/\Q)}$ such that ${\eta\vert_{K_{xy}}=\id}$  we have 
		${\gamma^{-1}\eta\gamma\vert_L=\eta}$. 
		
		We denote ${\eta^\gamma=\gamma^{-1}\eta\gamma}$. Proposition~\ref{pdihedral} implies that 
		$$
		\eta^\gamma\vert_{K_x(x)}\in \bigl\{\eta\vert_{K_x(x)}, \eta^{-1}\vert_{K_x(x)}\bigr\}. 
		$$
		We also have 
		${\eta^\gamma\vert_{K_y}= \id\vert_{K_y}=\eta\vert_{K_y}=\eta^{-1}\vert_{K_y}}$. 
		It follows that 
		$$
		\eta^\gamma\vert_{K_{xy}(x)}\in \bigl\{\eta\vert_{K_{xy}(x)}, \eta^{-1}\vert_{K_{xy}(x)}\bigr\}. 
		$$
		In particular, 
		${\eta^\gamma\vert_L\in \bigl\{\eta\vert_L, \eta^{-1}\vert_L\bigr\}}$. 
		Since ${\eta\vert_{K_{xy}}=\id}$ and ${L/K_{xy}}$ is $2$-elementary, we have ${\eta\vert_L= \eta^{-1}\vert_L}$. Hence ${\eta^\gamma\vert_L=\eta\vert_L}$. The proposition is proved. 
	\end{proof}

\subsubsection{(Almost) $2$-elementary  discriminants}
	
We will need a slight generalization of the notion of a $2$-elementary group. A finite abelian group~$G$ will be called  \textit{almost $2$-elementary} if it has a $2$-elementary subgroup of index~$2$. This means that either~$G$ is $2$-elementary, or it is ${\calC_4}$ times a $2$-elementary group. Recall that we denote by $\calC_m$ the cyclic group of order~$m$. 	
	
We call a discriminant (almost) $2$-elementary if  its class group is so.  (Almost) $2$-elementary discriminants can be conveniently characterized in terms of the $2$-rank, see Subsection~\ref{ssstwor}:
\begin{align}
\label{etwoel}
\text{$\Delta$ is $2$-elementary} &\Longleftrightarrow  h(\Delta)=2^{\rho_2(\Delta)};\\
\label{ealmtwoel}
\text{$\Delta$ is almost $2$-elementary} &\Longleftrightarrow  h(\Delta)\in \{2^{\rho_2(\Delta)},  2^{\rho_2(\Delta)+1}\}. 
\end{align}

	\begin{proposition}\label{pconductors}
Let ${D\ne -3,-4}$ be a fundamental discriminant, and let~$f$ be such that  ${Df^2}$ is an almost $2$-elementary discriminant. Then
\begin{equation}
\label{efdivides}
f\mid 2^4\cdot 3\cdot5\cdot7\cdot17.
\end{equation}		
	\end{proposition}

\begin{proof}
Denote ${\Delta=Df^2}$. Since ${D\ne -3,-4}$,  the class number formula~\eqref{eclnfr}, together with  equations~\eqref{epsild} and~\eqref{efactor}, implies that 
${h(\Delta)=h(D)\Psi}$, 
where 
$$
\Psi=\Psi(f,D)=f\prod_{p\mid f}\left(1-\frac{(D/p)}{p}\right). 
$$
To start with, note that, when~$\Delta$ is almost $2$-elementary,~$\Psi$ is a power of~$2$ by~\eqref{ealmtwoel}. More precisely, we have the following:
\begin{equation}
\label{eboundpsi}
\Psi\mid
2^{\omega(f)+2}  
\end{equation}
(where $\omega(\cdot)$ stands for the number of distinct prime divisors).
Indeed, if~$\Delta$ is almost  $2$-elementary, then so is~$D$, and we have both 
$$
h(D)\in \{2^{\rho_2(D)}, 2^{\rho_2(D)+1}\} \quad\text{and}\quad h(\Delta)\in \{2^{\rho_2(\Delta)} , 2^{\rho_2(\Delta)+1}\}.
$$
Hence,
\begin{equation}
\label{erhoael}
\nu_2(\Psi)\le  \rho_2(\Delta)-\rho_2(D)+1.  
\end{equation}
Proposition~\ref{pgauss} implies that 
\begin{equation}
\label{eboundrhotwominusrhotwo}
\rho_2(\Delta)-\rho_2(D)\le \omega(\Delta)-\omega(D)+1 \le \omega(f)+1,   
\end{equation}
which, together with~\eqref{erhoael}, proves~\eqref{eboundpsi}.

The following immediate observation is crucial and will be systematically used throughout the proof:
\begin{equation*}
p\mid f \Longrightarrow p-(D/p)\mid \Psi, \qquad
p^2\mid f \Longrightarrow p\mid \Psi. 
\end{equation*}
This implies very strong constraints on  the prime divisors of~$f$. 

%\paragraph{Prime divisors of~$f$}
First of all,~$D$ and~$f$ cannot have a common prime divisor other than~$2$. Indeed, if~$p$ divides both~$D$ and~$f$ then ${(D/p)=0}$ and ${p-(D/p)=p\mid \Psi}$. Since~$\Psi$ is a power of~$2$, we must have ${p=2}$. 

Next, if ${p\mid f}$ then either ${p+1}$ or ${p-1}$ is a power of~$2$.  
Indeed, if~$p$ is odd then ${p\nmid D}$ (as we have just seen) which implies that ${p-(D/p)\in \{p-1,p+1\}}$.

Yet another observation: 
if ${p^2\mid f}$ then ${p=2}$. Indeed, in this case we again have ${p\mid \Psi}$.

Thus, 
\begin{equation}
\label{eformoff}
f=2^kp_1\cdots p_m, 
\end{equation}
where~$k$ and~$m$ are non-negative integers and ${p_1, \ldots, p_m}$ are distinct odd primes not dividing~$D$. We claim that
\begin{equation}
\label{eboundk}
k=\nu_2(f) \le 4. 
\end{equation}
Indeed, if ${k\ge 1}$ then ${\omega(f) =m+1}$, and 
$$
2^{\omega(f)+2}=2^{m+3}\ge \Psi \ge 2^{k-1}(p_1-1)\cdots (p_m-1) \ge 2^{m+k-1}, 
$$ 
which proves~\eqref{eboundk}. 

%\paragraph{Prime divisors revisited}
To complete the proof, we have to show that the only possible prime divisors of~$f$ are $2,3,5,7,17$. If ${w(f) =1}$ then ${\Psi \mid 8}$, and ${f=2^k}$ or ${f=p}$, an odd prime.  Since ${p-1}$ or ${p+1}$ divides~$\Psi$, this implies that ${p\le 7}$.

Now assume that ${\omega(f)\ge 2}$, and that~$f$ has a prime divisor ${p\ne 2,3,5,7,17}$. Then ${p\ge 31}$, because one of ${p\pm1}$ must be a power of~$2$. 
Writing ${f=2^kp_1\cdots p_m}$ with ${2<p_1<\cdots <p_m}$ and ${p_m\ge 31}$, we have 
$$
2^{m+3}\ge 2^{\omega(f)+2} \ge \Psi  \ge (p_1-(D/p_1)) \cdots (p_m-(D/p_m)). 
$$
Since ${p_m\ge 31}$ and ${p_m-(D/p_m)}$ is a power of~$2$, we have ${p_m-(D/p_m)\ge 32}$. 

If ${m\ge 2}$ then 
$$
(p_1-(D/p_1)) \cdots (p_m-(D/p_m)) \ge (3-1)(5-1)^{m-2}\cdot 32 = 2^{2m+2}. 
$$
We obtain that ${m+3\ge 2m+2}$, which is impossible for ${m\ge 2}$. 
It follows that  ${m=1}$, in which case ${\omega(f)=2}$ and ${f=2^k p}$, where ${k\ge 1}$ and ${p\ge 31}$. We obtain  
${16\ge \Psi \ge 2^{k-1}\cdot 32}$, 
a contradiction. The proof is complete. 
\end{proof}

One can do some case-by-case analysis and show that~$f$ satisfies a stronger (but also more complicated) condition than~\eqref{efdivides}. However,~\eqref{efdivides} is sufficient for our purposes.  

We also need a similar result for the fundamental discriminants~$-3$ and~$-4$. 

\begin{proposition}
\label{pcondsttf}
If ${\Delta=-4f^2}$ is an almost $2$-elementary discriminant, then 
\begin{equation}
\label{ed-four}
\text{${f\le 8}$ or ${f\in \{10,12,15,20\}}$}. 
\end{equation}
In particular, ${h(\Delta)\mid 8}$.

If ${\Delta=-3f^2}$ is an almost $2$-elementary discriminant, then 
\begin{equation}
\label{ed-three}
\text{${f\le 5}$ or ${f\in \{7,8,11,13,16\}}$}. 
\end{equation}
In particular, ${h(\Delta)\mid 8}$.
\end{proposition}

\begin{proof}
Assume first that ${\Delta=-4f^2}$ is almost $2$-elementary. 
The class number formula~\eqref{eclnfr}, together with  equations~\eqref{epsild} and~\eqref{efactor}, implies that ${h(\Delta)=\Psi/2}$, 
where 
$$
\Psi=\Psi(f,-4)=f\prod_{p\mid f}\left(1-\frac{(-4/p)}{p}\right). 
$$
As in the proof of Proposition~\ref{pconductors}, this~$\Psi$ must be a power of~$2$; more precisely, 
\begin{equation}
\label{epsib}
\Psi \mid 2^{\omega(f)+2}. 
\end{equation}
Indeed, if ${2\nmid f}$, then  ${\Delta\equiv12\bmod 16}$ and ${\rho_2(\Delta) \le \omega(\Delta)-1=\omega(f)}$, while when ${2\mid f}$, we have ${\rho_2(\Delta) \le \omega(\Delta)=\omega(f)}$. 
In both cases we obtain 
$$
\Psi/2=h(\Delta) \mid 2^{\omega(f)+1}, 
$$
which is~\eqref{epsib}.  

Let~$p$ be an odd prime divisor of~$f$. As in the proof of Proposition~\ref{pconductors}, we have ${p^2\nmid f}$ and ${p-(-4/p)\mid \Psi}$; in particular, ${p-(-4/p)}$  is a power of~$2$. 
As in the proof of Proposition~\ref{pconductors}, we write the prime factorization of~$f$ as  ${f=2^kp_1\cdots p_m}$, where ${2<p_1<\cdots < p_m}$. We claim that
\begin{equation}
\label{ebkpm}
k+m\le 3 \quad\text{and}\quad m\le 2.  
\end{equation} 
Indeed, if ${k\ge 1}$ then ${\omega(f)=m+1}$, and~\eqref{epsib} implies that 
$$
2^{m+3} \ge \Psi \ge 2^{k-1}(2-(-4/2)) (3-(-4/3))^m = 2^{k+2m}, 
$$
which proves~\eqref{ebkpm} in the case ${k\ge 1}$. Similarly, if ${k=0}$ then 
$
{2^{m+2}\ge \Psi\ge 2^{2m}}
$,
proving~\eqref{ebkpm} in this case  as well. 

In a similar fashion one proves that 
\begin{equation}
\label{epmseven}
p_m\le 7.
\end{equation} 
Indeed, if ${p_m>7}$ then ${p_m\ge 17}$, because one of ${p_m\pm1}$ must be a power of~$2$. If ${k\ge 1}$ then
$$
2^{m+3} \ge \Psi \ge 2^k \cdot 4^{m-1}\cdot 16  = 2^{k+2m+2}, 
$$
which is impossible; if ${k=0}$ then 
$$
2^{m+2}\ge \Psi \ge 4^{m-1}\cdot 16= 2^{2m+2},
$$
again impossible. This proves~\eqref{epmseven}. 

It follows from~\eqref{ebkpm} and~\eqref{epmseven} that there are  finitely many possible~$f$. Checking them all using a \textsf{PARI} script, we obtain~\eqref{ed-four}.

\bigskip

Now assume that ${\Delta=-3f^2}$ is almost $2$-elementary. In this case ${h(\Delta)=\Psi/3}$, 
where 
$$
\Psi=\Psi(f,-3)=f\prod_{p\mid f}\left(1-\frac{(-3/p)}{p}\right). 
$$
This time~$\Psi$ must be~$3$ times a power of~$2$. It follows again that for an odd prime~$p$  we have ${p^2\nmid f}$. This is clear when ${p\ne 3}$, and if ${9\mid f}$ then 
$$
3(3-(-3/3))=9\mid\Psi,
$$
a contradiction. 

For every ${p\mid f}$ the difference ${p-(-3/p)}$ must be either a power of~$2$, or~$3$ times a power of~$2$. We claim that${p-(-3/p)}$ cannot be a power of~$3$. This is clear for ${p=2}$ and for ${p=3}$. Now assume that ${p\ne 2,3}$ and   ${p-(-3/p)=2^n}$. If~$n$ is odd then ${3\mid 2^n+1}$, which means that ${p=2^n-1}$. But in this case ${p\equiv 1\bmod 3}$ and ${p\equiv -1\bmod 4}$, which implies that ${(-3/p)=1}$. It follows that ${p-(-3/p)=2^n-2}$, a contradiction. Similarly, when~$n$ is even, we have ${p=2^n+1}$ and ${(-3/p)=-1}$, again a contradiction.

Thus, for every ${p\mid f}$ we have ${p-(-3/p)=3\cdot 2^n}$ for some~$n$. This implies that~$f$ cannot have two distinct prime divisors: if it did, then~$\Psi$ would be divisible by~$9$, a contradiction.

Thus, either ${f=2^k}$ for some~$k$, or ${f=p}$, an odd prime. This implies that ${\rho_2(\Delta) \le \omega(\Delta)\le 2}$, and 
$$
\Psi/3=h(\Delta) \mid 2^{\rho_2(\Delta)+1}\mid 8. 
$$
It follows that either ${f\mid 16}$, or ${f\in \{3,5,7,11,13,23\}}$. Checking all possible~$f$ using a \textsf{PARI} script,
we obtain~\eqref{ed-three}. 
\end{proof}

{\sloppy	

	\begin{proposition}
		\label{pbilly}
		There exists a fundamental discriminant~$D^\ast$ such that  ${h(D^\ast)\ge 128}$ and  the following holds. Let ${\Delta=Df^2}$ be either $2$-elementary or almost $2$-elementary. Then either ${D=D^\ast}$ or 
		$$
		h(\Delta) \le 
		\begin{cases}
		16, & \text{if $\Delta$ is $2$-elementary},\\
		64, & \text{if $\Delta$ is almost $2$-elementary}. 
		\end{cases}
		$$
	\end{proposition} 
	
}	
	
	This is proved in~\cite{ABP14}, an early version of~\cite{ABP15}, see Corollary~2.5 and Remark~2.6 therein.  This result was not included in the published version of~\cite{ABP15}, so we reproduce here the proof (adding some details missing in~\cite{ABP14}). The proof broadly follows the strategy in Weinberger \cite{We73}; in particular, it rests on a classical bound of Tatuzawa \cite[Theorem~2]{Ta51}. Given a fundamental discriminant~$D$, the \textit{$L$-function attached to~$D$} is $L(s, \chi)$, where~$\chi$ is the  quadratic character defined by the Kronecker symbol: ${\chi(n)=(D/n)}$. 
	
	\begin{lemma}[Tatuzawa]
	\label{ltatuzawa}
		Let $0<\varepsilon<1/2$. Then there exists a fundamental discriminant~$D^\ast$ such that the following holds. Let~$D$ be a fundamental discriminant and  $L(s,\chi)$  the attached $L$-function. Then ${ L(1,\chi) \geq 0.655 \varepsilon |D|^{-\varepsilon}}$ when ${|D| \leq \max\{e^{1/\varepsilon},73130 \}}$ and ${D\ne D^\ast}$. 
	\end{lemma}
	
	\begin{proof}[Proof of Proposition \ref{pbilly}]
If ${D=-3}$ or~$-4$ then the result  follows from Proposition~\ref{pcondsttf}. Hence we may assume that ${D\ne -3,-4}$. In this case  	the analytic class number formula states that ${h(D)=\pi^{-1}|D|^{1/2} L(1,\chi) }$.  		If $\Delta$ is almost $2$-elementary, then so is $D$. By~\eqref{ealmtwoel} and Proposition \ref{pgauss} we have ${h(D) \le 2^{\rho_2(D)+1}\le 2^{\omega(D)}}$.

		We pick $\varepsilon=0.048$ throughout and we use the corresponding~$D^\ast$ from Lemma~\ref{ltatuzawa}.  
Assuming that ${D\ne D^\ast}$, Lemma \ref{ltatuzawa} implies that
		$$ 
		 2^{\omega(D)}\geq h(D)=\pi^{-1}|D|^{1/2}L(1,\chi) \geq 0.655 \pi^{-1}\varepsilon |D|^{1/2-\varepsilon} 
		$$ 
as long as ${|D| \geq 1.2\cdot 10^9}$. This implies that 
		$$ |D| \leq \left((0.655\varepsilon)^{-1} \pi 2^{\omega(D)}  \right)^{1/(1/2-\varepsilon)}\leq 26549 \cdot 4.635^{\omega(D)};$$ we conclude that 
		\begin{equation}\label{econd1}
		|D| \leq \max\{ 26549 \cdot 4.635^{\omega(D)}, 1.2\cdot 10^9\}.
		\end{equation}
		Note moreover that, since $D$ is fundamental and 
$$
1.2\cdot 10^9<4\cdot(3 \cdot 5 \cdot 7\cdot 11  \cdots 37 )
$$ 
($4$ times the product of the first $11$ odd primes), we must have ${\omega(D) \leq 11}$ whenever ${|D| \leq 1.2\cdot 10^9}$. More generally, $|D|$ is at least~$4$ times the product of the first $\omega(D)-1$ odd primes. Hence, when  ${\omega (D) \geq 12}$, we have 
$$
|D|\ge 4\cdot(3 \cdot 5 \cdot 7\cdot 11  \cdots 37 ) \cdot 41^{\omega(D)-12}. 
$$
Combining this observation with the upper bound~\eqref{econd1}, we conclude that, when ${\omega(D) \geq 12}$,  we must have
		$$ 
		4\cdot(3 \cdot 5 \cdot 7\cdot 11 \cdots  37 ) \cdot 41^{\omega(D)-12} \leq |D| \leq 26549 \cdot 4.635^{\omega(D)}.
		$$
This is easily seen to be a contradiction for ${\omega(D) \geq 12}$. We conclude that  
\begin{equation}
\label{econd2} 
\omega(D) \leq 11
\end{equation} 
for any almost $2$-elementary fundamental ${D\ne D^\ast}$.
		
Thus, we are now left with the task of examining  discriminants ${\Delta=Df^2}$ such that  the corresponding fundamental~$D$ satisfies conditions~\eqref{econd1} and~\eqref{econd2}. We want to show that
\begin{itemize}
\item
if such~$\Delta$ is $2$-elementary then 	${h(\Delta)\le 16}$;

\item
if such~$\Delta$ is almost $2$-elementary then 	${h(\Delta)\le 64}$.
	
\end{itemize}

Proving this is a numerical check using \textsf{PARI}. We distinguish two cases: ${\omega(D)\le 6}$ and ${7\le\omega(D)\le 11}$. 

When ${\omega(D) \le 6}$ and~$D$ is almost $2$-elementary then ${h(D)\mid 64}$. Table~4 of Watkins~\cite{Wa04} implies that in this case  ${|D|\le 693067}$. Using  Proposition~\ref{pconductors} and a \textsf{PARI} script, we computed all $2$-elementary and all almost $2$-elementary discriminants ${\Delta=Df^2}$ such that ${|D|\le 693067}$. 
Our script 	found 101 discriminants that are $2$-elementary, the largest being ${-7392=-1848\cdot2^2}$. The class numbers of all these discriminants do not exceed $16$. Similarly, the script found  425 almost $2$-elementary discriminants,  ${87360=-5460\cdot4^2}$ being the largest, and their class numbers do not exceed $64$. % The total running time on a personal computer was less than 2 minutes. 
This completes the proof in the case ${\omega(D) \le 6}$.

When ${7\le\omega(D)\le 11}$, we can no longer use~\cite{Wa04}. %, and the only upper bound that we have is~\eqref{econd1}. %Checking all fundamental discriminants below this bound is too costly, so we proceed differently. 
To complete the proof in this case, for every ${n=7,\ldots,11}$ we determine all fundamental discriminants~$D$ satisfying 
\begin{equation}
\label{enappears}
\omega(D)=n, \qquad |D|\le 26549 \cdot 4.635^n
\end{equation}
(note that ${26549 \cdot 4.635^n>1.2\cdot10^9}$ for ${n\ge 7}$), and for each of them we check whether it is almost $2$-elementary. Our script found no  almost $2$-elementary fundamental discriminants satisfying~\eqref{enappears} with ${7\le n\le 11}$.   %The total running time was about 6 minutes for ${n=7}$, about 2 minutes for each of ${n=8,9}$, and negligible for ${n=10,11}$. 
This completes the proof of Proposition~\ref{pbilly}. 
	\end{proof}
	
	\begin{corollary}
		\label{cxy}
		Let~$x$ and~$y$ be singular moduli of distinct fundamental discriminants: ${D_x\ne D_y}$. 
		\begin{enumerate}
			\item
			\label{isame}
			Assume that ${\Q(x)= \Q(y)}$. Then  ${h(\Delta_x)= h(\Delta_y) \le 16}$. 
			
			\item
			\label{itwo}
			Assume that ${\Q(x)\subset \Q(y)}$ and ${[\Q(y):\Q(x)]=2}$. Then ${h(\Delta_x) \le 16}$ and ${h(\Delta_y) \le 32}$. 
		\end{enumerate}
	\end{corollary}
	
	\begin{proof}
		If ${\Q(x)= \Q(y)}$ then both $\Delta_x$ and~$\Delta_y$ are $2$-elementary by Proposition~\ref{pintersect}. Since ${D_x\ne D_y}$, one of the two is distinct from $D^\ast$; say, ${D_x\ne D^\ast}$. Then ${h(\Delta_x)\le 16}$. Hence ${h(\Delta_y) =h(\Delta_x) \le 16}$ as well. This proves item~\ref{isame}. 
		
		Now assume that we are in the situation of item~\ref{itwo}. Then $\Gal(\Q(x)/\Q)$  is $2$-elementary by Proposition~\ref{pintersect}. Hence so is $\Gal(K_y(x)/K_y)$. Since 
		$$
		[K_y(y):K_y(x)]\le [\Q(y):\Q(x)]=2,
		$$
the group $\Gal(K_y(y)/K_y)$ is  almost $2$-elementary. If ${D_x\ne D^\ast}$ then ${h(\Delta_x)\le 16}$ and ${h(\Delta_y) \le 32}$, so we are done. If ${D_y\ne D^\ast}$ then ${h(\Delta_y) \le 64}$. It follows that ${h(\Delta_x)\le 32}$ and we must have ${D_x\ne D^\ast }$, so we are done again. 
	\end{proof}

	\subsection{Gauss reduction theory, denominators}
	\label{ssgauss}

	Denote~$T_\Delta$ the set of triples  ${(a,b,c)\in \Z^3}$ with ${\Delta=b^2-4ac}$ satisfying
\begin{align}
	&
\text{$\gcd(a,b,c)=1$},\nonumber\\
	\label{ekuh}
	&
	\text{either\quad $-a < b \le a < c$\quad or\quad $0 \le b \le a = c$}. 
\end{align}
	Note that condition~\eqref{ekuh} is equivalent to
	$$
	\frac{b+\sqrt\Delta}{2a}\in \calF. 
	$$
	For every singular modulus~$x$ of discriminant~$\Delta$ there exists a unique triple ${(a_x,b_x,c_x)\in T_\Delta}$ such that, denoting 
	$$
	\tau_x= \frac{b_x+\sqrt\Delta}{2a_x},  
	$$
	we have ${x=j(\tau_x)}$. This is, essentially, due to Gauss; see \cite[Section~2.2]{BLP16} for more details.

	We will call~$a_x$ the \textit{denominator} of the singular modulus~$x$. 
	
	Note that, alternatively,~$\tau_x$ can be defined as the unique ${\tau\in \calF}$ such that ${j(\tau)=x}$.
	
	We will say that a positive integer~$a$ is a denominator for~$\Delta$ if it is a denominator of some singular modulus of discriminant~$\Delta$; equivalently, there exist ${b,c\in \Z}$ such that ${(a,b,c)\in T_\Delta}$.

	It will often be more convenient to use the notation $a(x)$, $b(x)$, $\tau(x)$ etc. instead   of~$a_x$,~$b_x$,~$\tau_x$, etc.   
	
	\begin{remark}
		\label{rparity}
		It is useful to note that~$b_x$ and~$\Delta_x$ are of the same parity: ${b_x\equiv \Delta_x\bmod 2}$. This is because ${\Delta_x=b_x^2-4a_xc_x\equiv b_x^2\bmod 4}$. 
	\end{remark}

	For every~$\Delta$ there exists exactly one singular modulus of discriminant~$\Delta$ and of  denominator~$1$, which will be called the \textit{dominant} singular modulus of discriminant~$\Delta$. Singular moduli with denominator~$2$ will be called \textit{subdominant}.

	\begin{proposition}
		\label{pcountden}
		Let~$\Delta$ be a discriminant. Then for every ${a\in\{2,3,4,5\}}$ there exist at most~$2$ singular moduli~$x$ with ${\Delta_x=\Delta}$  and ${a_x=a}$. For every  ${A\in \{13,18,30\}}$ there exists at most $S(A)$ singular moduli~$x$ with ${\Delta_x=\Delta}$ and ${a_x<A}$, where $S(A)$ is given in the following table:
		$$
		\begin{array}{r|ccc}
		A&13&18&30\\
		\hline
		S(A)&32&48&99
		\end{array}. 
		$$
	\end{proposition}
	
	\begin{proof}
		Let~$a$ be a positive integer. For a residue class ${r\bmod 4a}$ denote $B(r)$ the number of ${b\in \Z}$ satisfying 
		${-a< b\le a}$ and  ${b^2\equiv r\bmod 4a}$. Denote $s(a)$ the biggest of all $B(r)$:
		$$
		s(a)=\max\{B(r): r\bmod4a \}. 
		$$
		The number of triples ${(a,b,c)\in T_\Delta}$ with given~$a$ does not exceed  $B(\Delta)$; hence it  does not exceed $s(a)$ either. 
		A quick calculation shows that ${s(a)=2}$ for ${a\in \{2,3,4,5\}}$, and 
		$$
		\sum_{a<A}s(a)=S(A)
		$$
		for ${A\in \{13,18,30\}}$. The proposition is proved. 
	\end{proof}
	
	We will also need miscellaneous facts about the (non)-existence of singular moduli of some specific shape.
	The following proposition will be used in this article only for ${p=3}$. We, however, state it for general~$p$, for the sake of further applications.
	
	\begin{proposition}
		\label{poddpden}
		Let~$\Delta$ be a discriminant and~$p$ an odd prime number. 
		
		\begin{enumerate}
			\item
			\label{ionemodthree}
			Assume that ${(\Delta/p)=1}$. If ${|\Delta|\ge 4p^2-1}$  then~$\Delta$ admits exactly~$2$ singular moduli with denominator~$p$. More generally, if ${|\Delta|\ge 4p^k-1}$ then~$\Delta$ admits exactly~$2$ singular moduli with denominator~$p^k$. 
			
			\item
			\label{idivbyp}
			Assume that ${p^2\mid \Delta}$, and let~$a$ be a denominator for~$\Delta$. Then either ${p\nmid a}$ or ${p^2\mid a}$. In particular,~$p$ is not a denominator for~$\Delta$. 
		\end{enumerate}
	\end{proposition}
	
	\begin{proof}
		By Hensel's lemma, the assumption ${(\Delta/p)=1}$ implies that the congruence ${b^2\equiv \Delta\bmod p^k}$ has exactly two solutions~$b$ satisfying  ${0<b<p^k}$, and exactly one of these solutions satisfies ${b^2\equiv \Delta\bmod 4p^k}$. If~$b$ is this solution and ${|\Delta|\ge 4p^{k}-1}$ then the  two triples ${(p^k,\pm b, (b^2-\Delta)/4p^k)}$ belong to $T_\Delta$. This proves  item~\ref{ionemodthree}. 
		
		If ${p^2\mid \Delta}$ and ${p\mid a}$ then ${p\mid b}$ and ${p\nmid c}$. Hence ${p^2\mid 4ac=b^2-\Delta}$, which implies that ${p^2\mid a}$. This proves item~\ref{idivbyp}. 
	\end{proof}
	
	Here is an analogue of Proposition~\ref{poddpden} for ${p=2}$. 
	
	\begin{proposition}
		Let~$\Delta$ be a discriminant. 
		\label{pevenden}
		\begin{enumerate}

			\item
			\label{ionemodeight}
			Assume that ${\Delta\equiv 1\bmod 8}$. If ${|\Delta|>15}$, then~$\Delta$ admits exactly~$2$ subdominant singular moduli, which are ${j\bigl((\pm1+\sqrt\Delta)/4\bigr)}$. More generally, if ${|\Delta|\ge 4^{k+1}-1}$ then~$\Delta$ admits exactly~$2$ singular moduli with denominator~$2^k$.

			\item
			\label{inonemodeight}
			If ${\Delta\ne 1\bmod 8}$ then it admits at most one subdominant singular modulus.

			\item
			\label{i432hensel}
			Let~$\Delta$ satisfy ${\Delta\equiv 4\bmod 32}$ and ${|\Delta|\ge 252}$.  Then it admits exactly~$2$ singular moduli of denominator~$8$. These are  
			\begin{equation}
			\label{ebprime}
			j\left(\frac{\pm b'+\sqrt{\Delta/4}}8\right),
			\qquad 
			b'=
			\begin{cases}
			1, & \Delta\equiv 36\bmod 64,\\ 
			3, & \Delta\equiv 4\bmod 64.
			\end{cases}
			\end{equation}
			More generally, if ${k\ge 3}$ and ${|\Delta|\ge 4^{k+1}-4}$ then~$\Delta$ admits exactly~$2$ singular moduli with denominator~$2^k$. 
			
			\item
			\label{i432none}
			Let~$\Delta$ satisfy ${\Delta\equiv 4\bmod 32}$ and let~$a$ be a denominator for~$\Delta$. Then either~$a$ is odd, or ${8\mid a}$. In particular,~$2$,~$4$ and~$6$ are not denominators  for~$\Delta$.

			\item
			\label{i16none}
			Let~$\Delta$ be divisible by~$16$  and let~$a$ be a denominator for~$\Delta$. Then either~$a$ is odd, or ${4\mid a}$. In particular,~$2$ is not a denominator  for~$\Delta$.
			
			Furthermore,~$\Delta$ admits at most one singular modulus  with denominator~$4$. 
			
			\item
			\label{ieveryn}
			Assume that~$\Delta$ is even, but ${\Delta\ne 4\bmod32}$, and that ${|\Delta|>76}$. Then~$2$ or~$4$ is a denominator for~$\Delta$.

		\end{enumerate}

	\end{proposition}
	
	\begin{proof}
		Items~\ref{ionemodeight} and~\ref{i432hensel} are proved using Hensel's Lemma exactly like item~\ref{ionemodthree} of Proposition~\ref{poddpden}; we omit the details. 
		
		Item~\ref{inonemodeight} follows from \cite[Proposition~2.6]{BLP16}, and item~\ref{ieveryn} is \cite[Proposition~3.1.4]{BLP20}.
		Note that in~\cite{BLP20} denominators are called \textit{suitable integers}.

		We are left with items~\ref{i432none} and~\ref{i16none}.  If ${\Delta\equiv 4\bmod 32}$ and  ${(a,b,c)\in T_{\Delta}}$ with ${2\mid a}$  then ${2\|b}$ and~$c$ is odd. Hence
		${b^2\equiv 4\bmod 32}$, which implies that ${4ac\equiv 0\bmod 32}$. This shows that ${8\mid a}$, which proves item~\ref{i432none}. 
		
		Finally,  if ${16\mid \Delta}$ and ${2\|a}$ then ${2\mid b}$ and ${2\nmid c}$, which implies that 
		$$
		b^2=\Delta+4ac\equiv 8 \bmod 16, 
		$$
		a contradiction. This proves the first statement in item~\ref{i16none}. Similarly, if ${a=4}$ then  ${4\mid b}$ and ${2\nmid c}$; in particular, ${4ac\equiv 16 \bmod 32}$. Hence
		$$
		b=
		\begin{cases}
		0,& \text{if $\Delta\equiv 16 \bmod 32$},\\
		4,& \text{if $\Delta\equiv 0 \bmod 32$}. 
		\end{cases}
		$$
		Thus, in any case there is only one choice for~$b$, which proves the second statement in  item~\ref{i16none}. 
		The proposition is proved. 
	\end{proof}

	It is useful to note that the dominant singular modulus is real; in particular,  there exists at least one real singular modulus of every discriminant. This has the following consequence.

	\begin{proposition}
		\label{preal}
		Let~$x$ be a singular modulus, and let~$K$ be a subfield of $\Q(x)$. Assume that~$K$ is Galois over~$\Q$. Then~$K$ is a totally real field.  
	\end{proposition}
	
	Since $\Q(x)$ is  Galois over~$\Q$ when $\Delta_x$ is $2$-elementary, this implies that singular moduli of $2$-elementary discriminants are all real.

	\subsection{Isogenies}

	Let~$\Lambda$ and~$\Mu$ be lattices in~$\C$. We say that~$\Lambda$ and $\Mu$ are isogenous	if~$\Lambda$ is isomorphic to a sublattice of~$\Mu$.  
	
	More specifically, given a positive integer~$n$, 
	we say that~$\Lambda$ and~$\Mu$ are $n$-isogenous if~$\Mu$ has a sublattice~$\Lambda'$, isomorphic to~$\Lambda$, such that the quotient group ${\Mu/\Lambda'}$ is cyclic of order~$n$. This relation is symmetric. 
	
	Recall that we denote by~$\H$ the Poincaré plane and by~$\calF$ the standard fundamental domain. It is well-known that, for ${z,w\in \H}$, the lattices  ${\langle z,1\rangle}$ and ${\langle w,1\rangle}$ are $n$-isogenous if and only if there exists ${\gamma \in \Mu_2(\Z)}$ with coprime entries and determinant~$n$ such that ${w=\gamma(z)}$. Imposing upon~$z$ and~$w$ certain reasonable conditions, one may assume that the matrix~$\gamma$ is upper triangular. 
	
	\begin{proposition}
		\label{pisog}
		Let ${z, w \in \H}$ and let~$n$ be a positive integer. Assume that
		\begin{equation}
		\label{eimzgen}
		w\in \calF \quad\text{and}\quad \Im\, z \ge n. 
		\end{equation}
		Then the following two conditions are equivalent. 
		\begin{enumerate}
			\item
			\label{iisog}
			The lattices ${\langle z,1\rangle}$ and ${\langle w,1\rangle}$ are $n$-isogenous.
			\item
			\label{iabd}
			We have 
			$$
			w=\frac{pz+q}{s},
			$$
			where ${p,q,s\in \Z}$ satisfy 
			$$
			p,s>0, \qquad ps=n, \qquad \gcd(p,q,s)=1. 
			$$
		\end{enumerate}
	\end{proposition}
	
	\begin{proof}
		The implication \ref{iabd}~$\Rightarrow$~\ref{iisog} is trivial and does not require~\eqref{eimzgen}. 
		
		Now assume that condition~\ref{iisog} holds, and let ${\gamma\in M_2(\Z)}$ be a matrix with coprime entries and determinant~$n$ such that ${w=\gamma(z)}$. There exists ${\delta \in \SL_2(\Z)}$ such that 
		${\delta \gamma}$ is an upper triangular matrix. Replacing~$\delta$ by ${\bigl(\begin{smallmatrix}1&\nu\\0&1\end{smallmatrix}\bigr)\delta}$ with a suitable ${\nu\in \Z}$, we may assume that 
		${w':=\delta\gamma(z)}$ satisfies 
		\begin{equation}
		\label{erew}
		-\frac 12<\Re \, w'\le \frac12. 
		\end{equation}
		Write ${\delta\gamma= \bigl(\begin{smallmatrix}p&q\\0&s\end{smallmatrix}\bigr)}$. Replacing, if necessary,~$\delta$ by $-\delta$, we may assume that ${p,s>0}$.
		Since ${ps=n}$ and 
		$$
		w'=\frac{pz+q}{s},
		$$
		we only have to prove that ${w'=w}$.
		
		Note that
		$$
		\Im\, w' = \frac1s\Im\, z\ge 1 
		$$
		by~\eqref{eimzgen}. Together with~\eqref{erew} this implies that ${w'\in \calF}$.
		But~$w$ belongs to~$\calF$ as well, by our assumption~\eqref{eimzgen}. Since ${w'=\delta w }$ and each $\SL_2(\Z)$-orbit has exactly one point in~$\calF$, we must have ${w'=w}$. 
	\end{proof}
	
	We say that two singular moduli  are $n$-isogenous if,	writing ${x=j(\tau)}$ and ${y=j(\upsilon)}$, the lattices ${\langle\tau,1\rangle}$ and ${\langle\upsilon,1\rangle}$ are $n$-isogenous.
	
	Singular moduli~$x$ and~$y$ are $n$-isogenous if and only if ${\Phi_n(x,y)=0}$, where ${\Phi_n(X,Y)}$ denotes the modular polynomial of level~$n$. Since ${\Phi_n(X,Y)\in \Q[X,Y]}$,  being $n$-isogenous is preserved by Galois conjugation: for any ${\sigma\in \Gal(\bar\Q/\Q)}$ the singular moduli~$x^\sigma$ and~$y^\sigma$ are $n$-isogenous as long as~$x$ and~$y$ are.

	For a positive integer~$n$ define 
	$$
	\calQ(n)=\left\{\frac rs: r,s\in \Z, rs=n\right\}. 
	$$
	For example,
	$$
	\calQ(12)=\left\{\frac1{12},\frac13,\frac34,\frac43,3,12\right\}. 
	$$
	The following property is an immediate consequence of Proposition~\ref{pisog}. 
	
	\begin{corollary}
		\label{cisog}
		Let~$x$ and~$y$ be $n$-isogenous singular moduli. Assume that ${|\Delta_x|^{1/2}\ge 2na_x}$. Then ${(a_y/f_y)/(a_x/f_x)\in \calQ(n)}$. In particular, 
		$$
		\frac1n\le\frac{a_y/f_y}{a_x/f_x}\le n. 
		$$
		When ${n=p}$ is a prime number, we have 
		${(a_y/f_y)/(a_x/f_x)\in \left\{p,1/p\right\}}$. 
	\end{corollary}
	
	The following simple facts will be repeatedly used, often without special reference. 
	
	\begin{proposition}
		\label{pvarisog}
		Let~$x$ and~$y$ be singular moduli. 
		\begin{enumerate}
			\item
			\label{idenoms}
			Assume that ${\Delta_x=\Delta_y}$ and 
			${\gcd(a_x,a_y)=1}$.   
			Then~$x$ and~$y$ are ${a_xa_y}$-isogenous. 
			
			\item
			\label{ibothdom}
			Assume that 
			${a_x=a_y=1}$ and  ${\Delta_x/e_x^2=\Delta_y/e_y^2}$,
			where~$e_x$ and~$e_y$ are coprime positive integers. Then~$x$ and~$y$ are ${e_xe_y}$-isogenous.  
			
			\item
			\label{isubdomisog}
			Two subdominant singular moduli of the same discriminant are either equal or $4$-isogenous. 
		\end{enumerate}

	\end{proposition}
	
	\begin{proof}
		To prove item~\ref{idenoms}, note that 
		$$
		\tau_y=\frac{a_x\tau_x+(b_y-b_x)/2}{a_y}. 
		$$
		Since ${b_x\equiv b_y\bmod 2}$ (see Remark~\ref{rparity}), this proves that~$x$ and~$y$ are $a_xa_y$-isogenous.

		For item~\eqref{ibothdom} we have 
		$$
		\tau_x= \frac{b_x+e_x\sqrt\Delta}{2} \quad\text{and}\quad  \tau_y= \frac{b_y+e_y\sqrt\Delta}{2},
		$$
		where ${\Delta=\Delta_x/e_x^2=\Delta_y/e_y^2}$. 
		Hence 
		$$
		\tau_y=\frac{e_y\tau_x +(b_ye_x-b_xe_y)/2}{e_x}
		$$
		Remark~\ref{rparity} implies now that ${b_xe_y\equiv b_ye_x\bmod 2}$, and we conclude  that~$x$ and~$y$ are $e_xe_y$-isogenous.  
		
		To prove item~\ref{isubdomisog}, note that  distinct subdominant singular moduli of the same discriminant~$\Delta$ must be of the form ${j(\tau)}$ and ${j(\tau')}$, where 
		$$
		\tau=\frac{-1+\sqrt\Delta}{4} \quad\text{and}\quad \tau'=\frac{1+\sqrt\Delta}{4},
		$$
		see  Proposition~\ref{pevenden}:\ref{ionemodeight}. 
		We have 
		${\tau'= (2\tau+1)/2}$,
		which implies that ${j(\tau)}$ and ${j(\tau')}$ are $4$-isogenous. 
	\end{proof}
	
	\subsection{Galois-theoretic lemmas}
	
	In this subsection we collect some lemmas with Galois-theoretic flavor that will be repeatedly used in the proofs of Theorems~\ref{thpizthree} and~\ref{thprimel}.

	\begin{lemma}
		\label{lnorou}
		Let~$m$ be a positive integer and~$x$ a singular modulus. Then ${\Q(x)=\Q(x^m)}$. In other words: if~$x$ and~$y$ are distinct singular moduli of the same discriminant then ${x^m\ne y^m}$. 
	\end{lemma}
	
	\begin{proof}
		See \cite[Lemma~2.6]{Ri19}. 
	\end{proof}

	\begin{lemma}
		\label{lperm}
		Let~$x$ and~$y$ be distinct singular moduli of the same discriminant, ${K=K_x=K_y}$ their common CM field, ${L=K(x)=K(y)}$ the ring class field and ${\sigma \in \Gal(L/\Q)}$ a Galois morphism. Assume that~$\sigma$ permutes~$x$ and~$y$:
		$$
		x^\sigma=y \quad\text{and}\quad y^\sigma=x.
		$$
		Then~$\sigma$ is of order~$2$. 
	\end{lemma}
	
	\begin{proof}
		If ${\sigma \notin \Gal(L/K)}$ then it is of order~$2$ because every element of ${\Gal(L/\Q)}$ not belonging to ${\Gal(L/K)}$ is of order~$2$. And if ${\sigma \in \Gal(L/K)}$ then ${\sigma^2=1}$ because ${x^{\sigma^2}=x}$ and ${L=K(x)}$. 
	\end{proof}

	\begin{lemma}
		\label{lmainimproved}
		Let~$x$,~$y$ be distinct singular moduli of the same discriminant and let~$K$ and~$L$ be as in Lemma~\ref{lperm}. 
		Let~$F$ be a proper subfield of $\Q(x,y)$; we denote ${G=\Gal(L/F)}$. Then one of the following conditions is satisfied.
		
		\begin{enumerate}
			\item
			\label{ifixed}
			There exists ${\sigma\in G}$ such that, up to switching~$x,y$, we have ${x^\sigma=x}$ but ${y^\sigma\ne y}$. 
			
			\item
			\label{icyclictwo}
			We have ${[\Q(x,y):F]=2}$ and the non-trivial automorphism of $\Q(x,y)/F$ permutes~$x$ and~$y$.

			\item
			\label{icyclicthree}
			We have ${L=\Q(x,y)}$ and ${[L:F]=3}$. Moreover, there exists a singular modulus~$z$ and ${\sigma \in G}$ such that 
			$$
			x^\sigma=y, \quad y^\sigma=z, \quad z^\sigma=x. 
			$$
			
			\item
			\label{isigmaimproved}
			There exists ${\sigma \in G}$ such that 
			\begin{equation}
			\label{esigmaimproved}
			x^\sigma, y^\sigma\notin\{ x,y\} 
			\end{equation}
		\end{enumerate}
	\end{lemma}
	
	(Versions of this lemma were used, albeit implicitly, in~\cite{FR18} and  elsewhere, but  it does not seem to have appeared in the literature in this form.)

	\begin{proof}
		We may assume that every element in~$G$ which fixes~$x$ or~$y$ fixes both of them; otherwise we have item~\ref{ifixed}. We may also assume that~$x$ and~$y$ are conjugate over~$F$; otherwise, any ${\sigma \in  G}$ not belonging to ${\Gal(L/\Q(x,y))}$ satisfies~\eqref{esigmaimproved}.

		Assume first that ${L=\Q(x,y)}$. Then the only element of~$G$ that fixes~$x$ or~$y$ is the identity. Since~$x$ and~$y$ are  conjugate over~$F$, there is exactly one ${\sigma \in G}$ with the property ${x^\sigma= y}$ and exactly one ${\sigma'\in G}$ with the property ${y^{\sigma'}=x}$. Hence condition~\ref{isigmaimproved} holds if ${[L:F]\ge 4}$. And if ${[L:F]\le 3}$ then we have one of  conditions~\ref{icyclictwo} or~\ref{icyclicthree} is satisfied.  This completes the proof in the case ${L=\Q(x,y)}$.
		
		Now assume that ${L\ne\Q(x,y)}$. Since ${L=K(x)}$, the field $\Q(x)$ is a subfield of~$L$ of degree~$2$, and so is $\Q(y)$. If ${Q(x)\ne \Q(y)}$ then the compositum of these fields must be~$L$, which contradicts the assumption ${L\ne\Q(x,y)}$. Hence	
		 ${\Q(x)=\Q(y)}$ is a subfield of~$L$ of degree~$2$. 
		Condition~\ref{isigmaimproved} holds if ${[\Q(x):F]\ge 4}$, and condition~\ref{icyclictwo} does if ${[\Q(x):F]=2}$.
		
		We are left with the case ${[\Q(x):F]=3}$. In this case the Galois orbit of~$x$ over~$F$ consists of~$3$ elements: $x,y$ and a certain~$z$. 
		The group~$G$ must be either cyclic~$\calC_6$ or symmetric~$S_3$. In the latter case~$G$ acts by permutations on the set ${\{x,y,z\}}$. But in this case~$G$ has an element fixing~$x$ and permuting $y,z$, which is impossible because ${y\in \Q(x)}$.

		Thus, ${G=\calC_6}$. The group ${\Gal(L/\Q(x))}$ is a subgroup of~$G$; let~$\gamma$ be the non-trivial element  of ${\Gal(L/\Q(x))}$. Then ${\gamma \notin \Gal(L/K)}$; otherwise, from  ${L=K(x)}$ and ${x^\gamma=x}$ we would obtain that~$\gamma$ is the identity. It follows (see Proposition~\ref{pdihedral}) that~$\gamma$ does not commute with the elements of~$G$ of order~$3$, which is impossible, because~$G$ is an abelian group. The lemma is proved. 
	\end{proof}

	\begin{lemma}
		\label{ldeq}
		Let~$x$ and~$y$ be singular moduli with the same fundamental discriminant~$D$, and let ${K=\Q(\sqrt D)}$ be their common CM field. Assume that ${K(x)=K(y)}$, and that~$x$ and~$y$ do not both belong to~$\Q$.   Then  
		$$
		\Delta_x/\Delta_y\in \{4,1,1/4\}.
		$$ 
		Moreover, if (say) ${\Delta_x=4\Delta_y}$, then ${\Delta_y\equiv 1\bmod 8}$.
	\end{lemma}
	
	\begin{proof}
		See \cite[Proposition~4.3]{ABP15} and \cite[Subsection~3.2.2]{BLP16} (where the congruence ${\Delta_y\equiv 1\bmod 8}$ is proved). Note that in~\cite{ABP15} a formally stronger hypothesis ${\Q(x)=\Q(y)}$ is imposed, but in the proof it is only used that ${K(x)=K(y)}$. 
	\end{proof}
	
	\begin{lemma}
		\label{lndeq}
		Let $x,x',y,y'$ be singular moduli. Assume that
		$$
		\Delta_x=\Delta_{x'} \quad\text{and}\quad \Delta_y=\Delta_{y'}.
		$$
		Furthermore, assume that ${\Q(x,x')=\Q(y,y')}$. Then we have the following.
		\begin{enumerate}
			\item
			\label{inefudi}
			If ${D_x\ne D_y}$ then ${\Q(x)=\Q(y)}$.
			
			\item
			\label{iefudi}
			If ${D_x=D_y}$ then ${K(x)=K(y)}$, where ${K=K_x=K_y}$ is the common CM-field for~$x$ and~$y$.
		\end{enumerate}
	\end{lemma}
	
	\begin{proof}
		See \cite[Lemma~7.1]{BFZ20}. 
	\end{proof}

	\section{The linear relation}
	\label{slinear}
	
	Let ${x_1, \ldots, x_k}$ are non-zero singular moduli of the same fundamental discriminant~$D$ and ${m_1, \ldots, m_k\in \Z}$. We want to show that, under some reasonable assumption, the multiplicative relation 
	\begin{equation}
	\label{emultdep}
	x_1^{m_1}\cdots x_k^{m_k}=1
	\end{equation}
	implies the linear relation
	\begin{equation}
	\label{elinrelf}
	\sum_{i=1}^k\frac{f(x_i)}{a(x_i)}m_i=0.
	\end{equation}
	Recall that we denote~$f_x$ or $f(x)$ the conductor, and~$a_x$ or $a(x)$ the denominator of the singular modulus~$x$, as in Subsection~\ref{ssgauss}. 
	
Let us denote
	\begin{align}
	\label{edefxy}
	X=\max\{|\Delta(x_i)|: 1\le i\le k\} \quad\text{and}\quad
	Y=\min\{|\Delta(x_i)|: 1\le i\le k\}.
	\end{align}

	\begin{proposition}
		\label{plinearel}
		Let~$A$ be a positive number such that 
		\begin{equation}
		a(x_i)\le A \qquad (1\le i\le k). 
		\end{equation}
		Assume that 
		\begin{equation}
		\label{erootofy}
		Y^{1/2} > \frac13Ak(\log X+\log A+\log k+20). 
		\end{equation}
		Then~\eqref{emultdep} implies~\eqref{elinrelf}. 
	\end{proposition}

	It often happens that we control only a part of the denominators of ${x_1,\ldots, x_k}$. In this case we cannot expect an identity like~\eqref{elinrelf}, but we may have good bounds for the part of the sum corresponding to the terms with small denominators.    
	
	We need some extra notation. Set ${f=\gcd(f_{x_1}, \ldots, f_{x_k})}$ and ${\Delta=Df^2}$. We also define 
	\begin{align*}
	e_{x_i}=e(x_i)=f(x_i)/f \quad\text{and}\quad
	m_i'=e(x_i)m_i \qquad (i=1, \ldots, k).
	\end{align*}
	Then we have
	${\Delta(x_i)=e(x_i)^2\Delta}$, 
	and~\eqref{elinrelf} can be rewritten as 
	\begin{equation}
	\label{elinrele}
	\sum_{i=1}^k\frac{m_i'}{a(x_i)}=0.
	\end{equation}
	As indicated above, we want to obtain a less precise result, in the form of an inequality, which however holds true without the assumption that all the denominators are small. It will be practical to estimate separately the sums with positive and with negative exponents~$m_i$.

	\begin{proposition}
		\label{pineq}
		Let $A,\eps$ be real numbers satisfying ${A\ge 1}$  and  ${0<\eps \le 0.5}$. Assume that 
		\begin{equation}
		\label{eassumpdelta}
		|\Delta|^{1/2}\ge \max\left\{k\eps^{-1}\log X, \frac13A\bigl(\log(k\eps^{-1})+4\bigr)\right\}. 
		\end{equation}
		Then 
		\begin{align}
		\label{eupperforpos}
		\sum_{\overset{a(x_i)<A}{m_i>0}}\frac{m_i'}{a(x_i)}&\le \sum_{m_i<0} \frac{|m_i'|}{\min\{a(x_i),A\}} +\eps\|\bfm'\|, \\
		\label{eupperforneg}
		\sum_{\overset{a(x_i)<A}{m_i<0}}\frac{|m_i'|}{a(x_i)}&\le \sum_{m_i>0} \frac{m_i'}{\min\{a(x_i),A\}} +\eps\|\bfm'\|. 
		\end{align}
		Here we denote by $\|\bfm'\|$ the sup-norm  
		${\max\{|m_1'|, \ldots, |m_k'|\}}$. 
	\end{proposition}

	Propositions~\ref{plinearel} and~\ref{pineq} will be our principal tools in the proofs of Theorems~\ref{thpizthree} and~\ref{thprimel}. They will be
	proved in Subsection~\ref{sslinearproofs}, after some preparatory work in Subsections~\ref{ssest} and~\ref{ssmasser}.

	\subsection{Estimates for singular moduli}
	\label{ssest}

	For ${\tau\in \H}$ denote ${q=q_\tau:=e^{2\pi i\tau}}$. 
	Recall that the $j$-invariant function has the Fourier expansion
	$$
	j(\tau)=\sum_{k=-1}^\infty c_kq^k, 
	$$
	where 
	$$
	c_{-1}=1, \quad c_0=744, \quad c_1= 196884, \quad\ldots
	$$
	are positive integers. The fact that~$c_k$ are positive integers follows from the formula
$$
j(\tau) =q^{-1}\left(1+240\sum_{k=1}^\infty\sigma_3(k)q^k\right)^3\prod_{n=1}^\infty(1-q^n)^{-24}, \quad\text{where}\quad \sigma_3(k)=\sum_{d\mid k}d^3
$$	
(see, for instance, \cite[Chapter~1, Proposition~7.4 and Remark~7.4.1]{Si94}). 
Clearly, each of the series 
$$
\left(1+240\sum_{k=1}^\infty\sigma_3(k)q^k\right)^3 \quad\text{and}\quad (1-q^n)^{-24} \quad (n=1,2,3\ldots)
$$
has positive integer coefficients; hence so does the Fourier expansion of~$j$. 
		
	\begin{proposition}
\label{prestims}	
		Let ${\tau\in \H}$, and denote ${v=\Im\tau}$. 
\begin{enumerate}
\item
\label{ivgefive}
Assume that ${v\ge 5}$. Then 
\begin{align}
		j(\tau)
		\label{etwoterms}
		&=q^{-1}+744+ O_1(2\cdot10^5|q|), \\
		\label{ethreeterms}
		j(\tau)&=q^{-1}+744+196884q+ O_1(3\cdot10^7|q|^2),\\
		\label{elogsimple}
		\log |j(\tau)|&=2\pi v + O_1(800|q|)\\
		\label{elogoneterm}
		\log(qj(\tau))&=744q + O_1(5\cdot10^5|q|^2), \\
		\label{elogtwoterms}
		\log(qj(\tau))&=744q -79884q^2+O_1(2\cdot10^8|q|^3)
		\end{align}
\item
\label{ivlev}
Assume that ${\tau \in \calF}$ and ${v\le V}$, where~$V$ is a real number satisfying ${V\ge 5}$. 
Then 
		\begin{equation}
		\label{ecapv}
		\log |j(\tau)|\le 2\pi V + 3000e^{-2\pi V}. 
		\end{equation}
\end{enumerate}	
	\end{proposition}

We will use the following trivial lemma.

%	The following trivial lemma will be systematically used throughout the article. 

	\begin{lemma}
		\label{llogs}
		Let~$u$ be a complex number satisfying ${|u|<1}$. Then 
		\begin{equation*}
		|\log(1+u)|\le \frac{|u|}{1-|u|}.
		\end{equation*}
Furthermore, for ${n=1,2,\ldots}$ we have 
		\begin{equation*}
		\log(1+u)= \sum_{k=1}^n\frac{(-1)^{k-1}u^k}{k} +O_1\left(\frac{1}{n+1}\frac{|u|^{n+1}}{1-|u|}\right). 
		\end{equation*}
	\end{lemma}

	\begin{proof}[Proof of Proposition~\ref{prestims}] 
		Write ${\tau=u+vi}$. Then
$$
q=e^{2\pi ui}e^{-2\pi v} \quad\text{and}\quad
|q|= e^{-2\pi v}.
$$		
Let us prove item~\ref{ivgefive}. We have ${u\ge 5}$, which implies that
$$
|q|\le e^{-10\pi}. 
$$
For ${n\ge 0}$ denote
		$$
		j_n(\tau) = \sum_{k=n+1}^\infty c_kq^k. 
		$$
		In particular, 
		$$
		j_0(\tau) =j(\tau)-q^{-1}-744 \quad\text{and}\quad  j_1(\tau) =j(\tau)-q^{-1}-744-196884q. 
		$$
		Positivity of the coefficients~$c_k$ implies that  
		$$
		|j_n(\tau)q^{-n-1}|  \le\sum_{k=n+1}^\infty c_k|q|^{k-n-1}
		\le  \sum_{k=n+1}^\infty c_ke^{-10\pi (k-n-1)}= e^{10\pi(n+1)}j_n(5i). 
		$$
		In particular, 
		$$
		|j_0(\tau)q^{-1}|  \le e^{10\pi}j_0(5i)<2\cdot10^5 \quad\text{and}\quad |j_1(\tau)q^{-2}|  \le e^{20\pi}j_1(5i)<3\cdot10^7,
		$$
		which proves expansions~\eqref{etwoterms} and~\eqref{ethreeterms}. Using Lemma~\ref{llogs},   we deduce from them expansions~\eqref{elogoneterm} and~\eqref{elogtwoterms}. Furthermore,~\eqref{elogsimple} easily follows from~\eqref{elogoneterm}. 
		
Now let us prove item~\ref{ivlev}. We no longer have ${v\ge 5}$. However, since ${\tau \in \calF}$, we have		 ${v\ge \sqrt3/2}$. Hence 
		$$
		|j(\tau)|\le |q|^{-1}+ 744+j_0(\sqrt3/2) \le e^{2\pi v}+ 2079\le e^{2\pi V}(1+2079e^{-2\pi V}). 
		$$
		Using Lemma~\ref{llogs}, we obtain~\eqref{ecapv}. 
	\end{proof}
	
We want to apply Proposition~\ref{prestims} to singular moduli. If~$x$ is a singular modulus, then there exists a unique ${\tau_x\in \calF}$ such that ${x=j(\tau_x)}$. 	We have
$$
\tau_x= \frac{b+\sqrt\Delta}{2a}, 
$$
where ${\Delta=\Delta_x}$ is the discriminant,  ${a=a_x}$ is the denominator of the singular modulus~$x$, and ${b\in \Z}$; 	see Subsection~\ref{ssgauss} for more details. 
	
	\begin{corollary}
		\label{clogx}
		Let~$x$ be a singular modulus of discriminant~$\Delta$ and denominator~$a$.  
		Then we have the following. 
		
		\begin{enumerate}
			\item
		\label{iequ}
			Assume that  ${a\le 0.1|\Delta|^{1/2}}$. Then 
			$$
			\log|x|= \pi \frac{|\Delta|^{1/2}}{a}+O_1(e^{-2|\Delta|^{1/2}/a}). 
			$$		
			\item
			\label{inequ}
			Let ${A\ge1}$ be such that ${a\ge A}$ and ${A\le 0.1|\Delta|^{1/2}}$. Then
\begin{equation}
\label{elogxlev}
			\log|x|\le \pi \frac{|\Delta|^{1/2}}{A}+e^{-2|\Delta|^{1/2}/A}. 
\end{equation}
			In particular, if ${|\Delta|\ge 10^4}$ then 
\begin{equation}
\label{elogxleone}
			\log|x|\le \pi |\Delta|^{1/2}+e^{-2|\Delta|^{1/2}} \le 4|\Delta|^{1/2}. 
\end{equation}
\item
			\label{icomplog}
			Define~${\tau_x}$
			as in Subsection~\ref{ssgauss}.   
			Assume that ${a\le 0.1|\Delta|^{1/2}}$. Then 
			\begin{align}
			\label{efirst}
			\log(xq)& =744q + O_1(5\cdot10^5|q|^2), \\
			\label{esecond} 
			\log(xq)&=744q -79884q^2+O_1(2\cdot10^8|q|^3). 
			\end{align}
			where  ${q= e^{2\pi i\tau_x} }$.
		\end{enumerate}
\end{corollary}

\begin{proof}
Hypothesis ${a\le 0.1|\Delta|^{1/2}}$ implies that 
$$
\Im\tau_x = \frac{|\Delta|^{1/2}}{2a}\ge 5. 
$$
Applying~\eqref{elogsimple} with ${\tau=\tau_x}$, we obtain
$$
\log|x|= \pi \frac{|\Delta|^{1/2}}{a}+O_1(800e^{-\pi|\Delta|^{1/2}/a}). 
$$	
Since ${|\Delta|^{1/2}/a\ge 10}$, we have 
${800e^{-\pi|\Delta|^{1/2}/a} \le e^{-2|\Delta|^{1/2}/a}}$. This proves item~\ref{iequ}.

Similarly, item~\ref{inequ} follows by applying~\eqref{ecapv} with ${V= |\Delta|^{1/2}/2A}$. Note that~\eqref{elogxleone} is a special case of~\eqref{elogxlev}  corresponding to ${A=1}$.  

Finally,~\eqref{efirst} and~\eqref{esecond} are obtained by setting ${\tau=\tau_x}$ in~\eqref{elogoneterm} and~\eqref{elogtwoterms}, respectively. 
\end{proof}

	We also need a lower bound. The following  is (a weaker version of) \cite[Theorem~6.1]{BFZ20}, applied with ${y=0}$.  
	
	\begin{proposition}
		\label{pbifazhu}
		Let~$x$ be a singular modulus with discriminant ${\Delta_x\ne -3}$. Then ${|x|\ge|\Delta_x|^{-3}}$.
	\end{proposition}

	\subsection{Bounding the exponents}
	\label{ssmasser}

	Let ${\alpha_1, \ldots, \alpha_k}$ be non-zero algebraic numbers. The set of 
	${\bfm=(m_1, \ldots, m_k) \in \Z^k}$ such that 
	$$
	\alpha_1^{m_1}\cdots \alpha_k^{m_k}=1
	$$ 
	is a subgroup in $\Z^k$, denoted here by ${\Gamma(\alpha_1, \ldots,\alpha_k)}$ (or simply~$\Gamma$ if this does not cause confusion). 
	
	Masser~\cite{Ma88} showed that~$\Gamma$  admits a small $\Z$-basis. To state his result, let us introduce some notation. Let~$L$ be a number field. 
	We denote by ${\omega=\omega(L)}$ the order of the group of roots of unity belonging to~$L$, and by ${\eta=\eta(L)}$ the smallest positive height of the elements of~$L$:
	$$
	\eta=\min\{\height(\alpha): \alpha \in L,\ \height(\alpha)>0\}. 
	$$
	Here $\height(\cdot)$ is the usual absolute logarithmic height.\footnote{There is no risk of confusing the height $\height(\cdot)$ and the class number $h(\cdot)$, not only because the former is roman and the latter is italic, but, mainly,  because class numbers do not occur in this section, and heights do not occur outside this section.} 	
	Finally, we define the norm of a vector  ${\bfm=(m_1,\ldots,m_k)\in \Z^k}$  as
	${\|\bfm\|=\max\{|m_1|, \ldots, |m_k|\}}$. 
	
	\begin{proposition}[Masser]
		\label{pmasser}
		Let ${\alpha_1, \ldots\alpha_k}$ be elements in ${L^\times}$.
		Denote 
		$$
		\height=\max\{\height(\alpha_1), \ldots, \height(\alpha_k),\eta\}. 
		$$
		Then ${\Gamma(\alpha_1, \ldots,\alpha_k)}$ 
		has a $\Z$-basis consisting of vectors with norm bounded by 
		${\omega(k\height/\eta)^{k-1}}$. 
	\end{proposition}

	We want to adapt this result to the case when our algebraic numbers are singular moduli.
	
	\begin{proposition}
		\label{pbasrel}
		Let ${x_1, \ldots, x_k}$  be non-zero singular moduli. Set 
		$$
		X=\max\{|\Delta_{x_1}|,\ldots, |\Delta_{x_k}|\}, 
		$$
		and assume that, among ${D_{x_1},\ldots, D_{x_k}}$,  there are~$\ell$ distinct fundamental discriminants; in symbols:
$$
\ell=\#\{D_{x_1},\ldots, D_{x_k}\}. 
$$ 
Then the group $\Gamma(x_1, \ldots, x_k)$ has a $\Z$-basis consisting of vectors with norm bounded by 
		${ 24 (c(\ell)kX^{1/2})^{k-1}}$,
		where ${c(\ell)=3^{4^\ell+2^{\ell+1}+8}}$. In particular, %in the case when ${x_1, \ldots, x_k}$ all have the same fundamental discriminant, we have 
${c(1)=3^{16}}$. 
	\end{proposition}

	The proof uses the following result due to Amoroso and Zannier~\cite[Theorem~1.2]{AZ10}. 
	
	\begin{lemma}
		\label{laz}
		Let~$K$ be a number field of degree~$d$, and let~$\alpha$ be an algebraic number such $K(\alpha)$ is an abelian extension of~$K$. Then either ${\height(\alpha)=0}$ or ${\height(\alpha)\ge 3^{-d^2-2d-6}}$. 
	\end{lemma}
	
	\begin{proof}[Proof of Proposition~\ref{pbasrel}]
		Denote 
		$$
		K=\Q(\sqrt{D_{x_1}}, \ldots, \sqrt{D_{x_k}}), \qquad L=K(x_1, \ldots, x_k). 
		$$
		To apply Proposition~\ref{pmasser},
		we have to estimate the quantities~$\height$,~$\eta$ and~$\omega$.

		Since~$x_i$ is an algebraic integer, and every one of its conjugates~$x_i^\sigma$ satisfies  ${\log|x_i^\sigma|\le 4|\Delta|^{1/2}}$ (see  Corollary~\ref{clogx}:\ref{inequ}), 
		we have 
		${\height(x_i) \le 4X^{1/2}}$. 
		It follows that  
		${\height\le 4X^{1/2}}$. 
		
		Since ${[K:\Q]\le 2^\ell}$ and~$L$ is an abelian extension of~$K$, Lemma~\ref{laz} implies that  
		${\eta \ge 3^{-4^\ell-2^{\ell+1}-6}}$. Finally, we have ${\omega\le 24}$ from Corollary~\ref{crofoneinrcf}.
		
		Putting all of these together, the result follows. 
	\end{proof}

	In this article we will often work with relations of the form 
	$$
	x_1^{m_1}\cdots x_k^{m_k}= (x_1')^{m_1}\cdots (x_k')^{m_k}.
	$$

	It is useful to have a bounded basis for the group of these relations as well.
	
	\begin{proposition}
		\label{pbasrelsigma}
		Let ${x_1, \ldots, x_k, x_1', \ldots, x_k'}$ be singular moduli of discriminants not exceeding~$X$, and let~$\ell$ be the numbers of distinct fundamental discriminants among ${D_{x_1}, \ldots, D_{x_k'}}$.   
		Then the group $\Gamma(x_1/x_1', \ldots, x_k/x_k')$ has a $\Z$-basis consisting of vectors with norm bounded by 
		${ 24 (c(\ell)kX^{1/2})^{k-1}}$,
		where ${c(\ell)=3^{4^\ell+2^{\ell+1}+8}}$. In particular, ${c(1)=3^{16}}$. 
	\end{proposition}
	
	\begin{proof}
		Same as for Proposition~\ref{pbasrel}, only with ${\height\le 8X^{1/2}}$. 
	\end{proof}

	\subsection{Proofs of Propositions~\ref{plinearel} and~\ref{pineq}}
	
	\label{sslinearproofs}

	\begin{proof}[Proof of Proposition~\ref{plinearel}]
		By Proposition~\ref{pbasrel}  we may assume that
		\begin{equation}
		\label{embounded}
		\|\bfm\|\le  24 (3^{16}kX^{1/2})^{k-1}. 
		\end{equation}
		Using Corollary~\ref{clogx}, we
		obtain
		\begin{equation}
		\label{ezero=}
		0=\sum_{i=1}^km_i\log|x_i|= \pi|D|^{1/2}L+O_1(k\|\bfm\|e^{-3Y^{1/2}/A}), 
		\end{equation} 
		where~$L$ is the left-hand side of~\eqref{elinrelf}. (Recall that~$D$ denotes the common fundamental discriminant of ${x_1, \ldots, x_k}$.) Using~\eqref{erootofy} and~\eqref{embounded}, we deduce from this the estimate
		${|L|\le 0.5A^{-k}}$. Since~$L$ is a rational number with denominator not exceeding~$A^k$, we must have ${L=0}$. 
	\end{proof}

	\begin{proof}[Proof of Proposition~\ref{pineq}]
		We will prove only~\eqref{eupperforpos}, because~\eqref{eupperforneg} is totally analogous. 
		
		Using Corollary~\ref{clogx} and Proposition~\ref{pbifazhu}, we obtain
		\begin{align*}
		0&=\frac{1}{\pi|\Delta|^{1/2}}\sum_{i=1}^km_i\log |x_i| \\
		&\ge \sum_{\overset{1\le i\le k}{a(x_i)<A}}\frac{m_i'}{a(x_i)} - \sum_{m_i<0} \frac{|m_i'|}{\min\{a(x_i),A\}}+O_1\left(k\|\bfm'\| \frac{3\log X+e^{-3|\Delta|^{1/2}/A}}{\pi|\Delta|^{1/2}}\right),
		\end{align*}
		Using our hypothesis~\eqref{eassumpdelta}, we obtain
		\begin{align*}
		\frac{3k\log X}{\pi|\Delta|^{1/2}}&\le \frac3\pi\eps,\qquad \frac{ke^{-3|\Delta|^{1/2}/A}}{\pi|\Delta|^{1/2}} \le 0.01\eps, 
		\end{align*}
		and the result follows. 
	\end{proof}

	\section{Proof of Theorem~\ref{thprimel}}
	\label{sprimel}
	In this section we prove Theorem~\ref{thprimel}. Thus, throughout this section, unless the contrary is stated explicitly, 	$x$ and~$y$ are distinct singular moduli	and ${m,n}$ are non-zero integers such that  ${\Q(x^my^n)\ne \Q(x,y)}$. 
	
The proof is organized as follows. Assuming that	
	\begin{equation}
	\label{elowerprimelsix}
	\max\{|\Delta_x|,|\Delta_y|\}\ge 10^6,
	\end{equation}
we show that one of the following two conditions is satisfied: either
	\begin{equation}
	\label{ewewant}
	\begin{aligned}
	&\Delta_x=\Delta_y, \qquad m=n, \qquad [\Q(x,y):\Q(x^my^m)]= 2, \\
	&\text{$x$ and~$y$ are conjugate over $\Q(x^my^m)$} 
	\end{aligned}
	\end{equation}
(as wanted), or 
\begin{equation}
\label{edfd}
\text{${\{\Delta_x,\Delta_y\}=\{\Delta,4\Delta\}}$ for some ${\Delta\equiv 1\bmod 8}$}. 
\end{equation}
Unfortunately, we cannot rule out~\eqref{edfd} assuming merely~\eqref{elowerprimelsix}, but we show that~\eqref{edfd} is impossible under a stronger hypothesis 
\begin{equation}
	\label{elowerprimeleight}
	\max\{|\Delta_x|,|\Delta_y|\}\ge 10^8,  
  	\end{equation}	
completing thereby the proof.

	We assume that $x,y$ have the same fundamental discriminant; if this is not the case, then the argument is much simpler, see Subsection~\ref{ssdisfudis}.  We denote by~$K$ the common CM field of $x,y$, and we set 
	${L=K(x,y)}$. We also denote 
	$$
	\alpha=x^my^n, \qquad F=\Q(\alpha), \qquad G= \Gal (L/F).
	$$
	Since~$F$ is a proper subfield of $\Q(x,y)$, there exists ${\sigma \in G}$ such that ${x^\sigma\ne x}$ or ${y^\sigma \ne y}$. We claim that 
	\begin{equation}
	\label{ebothxy}
	x^\sigma \ne x \quad \text{and}\quad y^\sigma \ne y. 
	\end{equation}
	Indeed, if, say, ${y^\sigma=y}$ then ${(x^\sigma)^m=x^m}$, which implies ${x^\sigma=x}$ by Lemma~\ref{lnorou}. 
	
In the course of the argument we study study  multiplicative relations
	\begin{equation}
	\label{emainmultrelation}
	x^my^n(x^\sigma)^{-m}(y^\sigma)^{-n}=1,
	\end{equation}
	with various choices of ${\sigma \in G}$ satisfying~\eqref{ebothxy}, 
	using Propositions~\ref{plinearel} and~\ref{pineq}. 
In our usage of Propositions~\ref{plinearel} and~\ref{pineq}  the  parameters therein will satisfy the following restrictions:
	\begin{equation}
	\label{eparameters}
	\begin{split}
	k\le 4, \quad X=\max\{|\Delta_x|,|\Delta_y|\}\ge 10^6,\ \text{respectively},\ 10^8, \\ 
	Y\ge \frac14X, \quad A\le 9, \quad \eps=0.16,\ \text{respectively},\ 0.016. 
	\end{split}
	\end{equation}
	It is easy to verify that for any choice of parameters satisfying~\eqref{eparameters}, conditions~\eqref{erootofy} and~\eqref{eassumpdelta} are met, so using the propositions is justified.  
	
\begin{remark}
\label{repses}
Everywhere throughout the proof until Subsection~\ref{sssfourdems} we assume~\eqref{elowerprimelsix}, and we use Proposition~\ref{pineq} with ${\eps=0.16}$. Starting from Subsection~\ref{sssfourdems} we have~\eqref{edfd} and  assume~\eqref{elowerprimeleight}, which will allow us to use Proposition~\ref{pineq} with ${\eps=0.016}$. 
\end{remark}

	\subsection{A special case}
	\label{ssspecialcase}
	In this subsection we study the special case 
	\begin{equation}
	\label{especial}
	m=-n \quad\text{and}\quad\Delta_x=\Delta_y.
	\end{equation}
	We will need the result on this case to treat the general case. It will also be a good illustration of how our method works in a simple set-up. 
	
	Let ${\sigma \in G}$ be such that~\eqref{ebothxy} holds. We will apply Proposition~\ref{pineq} to the multiplicative relations 
	\begin{align}
	\label{erelwithsigma}
	x^my^{-m}(x^\sigma)^{-m}(y^\sigma)^m&=1,\\
	\label{erelwithsigmaminone}
	x^my^{-m}(x^{\sigma^{-1}})^{-m}(y^{\sigma^{-1}})^m&=1. 
	\end{align}
	We may assume, up to Galois conjugation, that $x$ is dominant. Then neither of $y,x^\sigma, x^{\sigma^{-1}}$ is:
	$$
	a(x)=1, \qquad a(y),\ a(x^\sigma),\ a(x^{\sigma^{-1}})\ge 2. 
	$$
	If one of ${a(y), a(x^\sigma)}$ is $\ge 3$, then, applying Proposition~\ref{pineq} to~\eqref{erelwithsigma} with ${A=3}$ and ${\eps=0.16}$, we obtain 
	$$
	m \le \left(\frac{1}{\min\{3, a(y)\}}+ \frac{1}{\min\{3, a(x^\sigma)\}} + 0.16\right)m \le \left(\frac12+\frac13+0.16\right)m, 
	$$
	a contradiction. 
	
	Thus, we must have ${a(y)=a(x^\sigma)=2}$. This means that~$y$ and~$x^\sigma$ are  either equal or $4$-isogenous, see  Proposition~\ref{pvarisog}:\ref{isubdomisog}. Hence so are $y^{\sigma^{-1}}$ and~$x$. Corollary~\ref{cisog} now implies  that 
	${a(y^{\sigma^{-1}})\le 4}$.

	Applying Proposition~\ref{pineq} to~\eqref{erelwithsigmaminone} with ${A=5}$ and ${\eps=0.16}$, we obtain 
	\begin{align*}
	m+\frac{m}{a(y^{\sigma^{-1}})} &\le m\left(\frac{1}{\min\{5,a(y)\}}+ \frac{1}{\min\{5,a(x^{\sigma^{-1}})\}}+0.16\right)\\
	& \le m\left(\frac12+\frac12+0.16\right). 
	\end{align*}
	Since ${a(y^{\sigma^{-1}})\le 4}$, it is a contradiction. 
	We have proved that~\eqref{especial} is impossible.

	\subsection{The general case: preparations}
	\label{ssprepage}
	Now we are ready to treat the general case. 
	Pick ${\sigma \in G}$ satisfying~\eqref{ebothxy}. 
	We have ${(x/x^\sigma)^m=(y^\sigma/y)^n}$, and, in particular, 
	$
	{\Q((x/x^\sigma)^m)=\Q((y/y^\sigma)^n)}
	$.
	The special case of  Theorem~\ref{thprimel}, treated in Subsection~\ref{ssspecialcase},  implies that 
	$$
	\Q((x/x^\sigma)^m)=\Q(x,x^\sigma), \qquad \Q((y/y^\sigma)^n)=\Q(y,y^\sigma). 
	$$
	Lemma~\ref{lndeq}:2 now implies that ${K(x)=K(y)=L}$, and Lemma~\ref{ldeq} implies that there exists a discriminant~$\Delta$ such that
	$$
	\Delta_x=e_x^2\Delta, \quad \Delta_y=e_y^2\Delta, \qquad (e_x,e_y)\in \{(1,1), (2,1),(1,2)\}, 
	$$
	and, moreover, 
	\begin{equation}
	\label{emoreover}
	\text{if $(e_x,e_y) \ne (1,1)$ then $\Delta\equiv 1\bmod 8$}. 
	\end{equation}
	We may and will assume in the sequel that 
	\begin{equation}
	\label{ebasicassum}
	m>0, \qquad e_xm\ge e_y|n|, \qquad a_x=1. 
	\end{equation}
	If ${(e_x,e_y) \ne (1,1)}$ then~$x$ and~$y$ are not conjugate over~$\Q$, and~\eqref{ebothxy} becomes
	\begin{equation}
	\label{emainoption}
	x^\sigma,y^\sigma \notin \{x,y\}. 
	\end{equation}
	When ${(e_x,e_y)=(1,1)}$,  Lemma~\ref{lmainimproved} implies that~$\sigma$ can be redefined to satisfy either~\eqref{emainoption} or one of the following: %  two more cases are possible:
	\begin{align}
	\label{etwocyc}
	&x^\sigma=y, \quad y^\sigma =x, \quad \hphantom{z^\sigma=x,} \qquad [\Q(x,y):F]=2;\\
	\label{ethreecyc}
	&x^\sigma=y, \quad y^\sigma=z,\quad z^\sigma=x, \qquad L=\Q(x,y), \quad [L:F]=3. 
	\end{align}

	Case~\eqref{etwocyc} is easy: relation~\eqref{emainmultrelation} becomes ${x^{m-n}=y^{m-n}}$, and Lemma~\ref{lnorou} implies that ${m=n}$, which means that we have~\eqref{ewewant}.
	
	We have to show that the other two cases are impossible. For~\eqref{ethreecyc} this is done in Subsection~\ref{sscycl}. Case~\eqref{emainoption} is much harder to dispose of, we deal with  it in Subsection~\ref{ssmainopone}.

	\subsection{Case~\eqref{ethreecyc}}
	\label{sscycl}
	
	We have
	$$
	x^my^{n-m}z^{-n}=1. 
	$$ 
	Recall that ${m\ge |n|}$ and  ${a_x =1}$, see~\eqref{ebasicassum}. In particular, 
	we must have ${a_y,a_z\ge 2}$. 
	
	Assume first that ${n>0}$.  Then
	$$
	\max\{m,|n-m|, |-n|\}=m. 
	$$
	Using Proposition~\ref{pineq} with ${\eps=0.16}$ and ${A=5}$, we obtain 
	$$
	m \le 
	\frac{m-n}{\min\{5,a_y\}}+\frac{m}{\min\{5,a_z\}}n+0.16m \le \frac12(m-n)+\frac12n+0.16m,
	$$
	a contradiction. 
	
	{\sloppy
		
		Now assume that ${n<0}$. Then
		$$
		\max\{m,|n-m|, |-n|\}\le 2m. 
		$$
		If ${a_y\ge 3}$ then Proposition~\ref{pineq} with ${\eps=0.16}$ and ${A=3}$ implies that 
		$$
		m\le \frac{1}{3}\cdot2m+0.16\cdot 2m, 
		$$
		a contradiction. If ${a_y=2}$ then~$x$ and~$y$ are  $2$-isogenous (see  Proposition~\ref{pvarisog}:\ref{idenoms}), and so are ${y=x^\sigma}$ and ${z=y^\sigma}$. This implies that ${a_z\in \{1,4\}}$, see Corollary~\ref{cisog}. But ${a_z\ge 2}$, and so ${a_z=4}$. Proposition~\ref{plinearel} now implies that 
		$$
		m+ \frac{n-m}{2}-\frac{n}{4}=0,
		$$
		yielding ${n=-2m}$, again a contradiction. Thus,~\eqref{ethreecyc} is impossible.  
		
	}

	\subsection{Case~\eqref{emainoption}}
	\label{ssmainopone}
	We will use notation 
	$$
	m'=e_xm, \qquad n'=e_yn. 
	$$
	Recall that ${m'\ge |n'|}$ and ${a(x)=1}$, see~\eqref{ebasicassum}. This implies, in particular, that ${a(x^\sigma)\ge 2}$.

	\subsubsection{One of~$y$,~$y^\sigma$ is dominant}

	We start by showing that either~$y$ or~$y^\sigma$ is dominant. 
	
	\begin{proposition}
		\label{pinside}
		If ${n>0}$ then ${a(y^\sigma)=1}$ and ${\sigma^2=1}$. If ${n<0}$ then ${a(y)=1}$. In both cases we have ${(e_x,e_y)\ne (1,1)}$ and ${\Delta\equiv 1\bmod 8}$.
	\end{proposition}
	
	\begin{proof}
		We treat separately  ${n>0}$ and ${n<0}$.

		\bigskip

		Assume first that ${n>0}$, but ${a(y^\sigma) \ge 2}$. 
		We know already that ${a(x^\sigma) \ge 2}$. If one of $a(x^\sigma),a(y^\sigma)$ is $\ge 3$ then, applying Proposition~\ref{pineq} with ${A=3}$ and ${\eps=0.16}$  to the relation
		${x^my^n(x^\sigma)^{-m}(y^\sigma)^{-n}=1}$, 
		we obtain
		$$
		m'\le \frac{m}{\min\{3,a(x^\sigma)\}}+ \frac{n}{\min\{3,a(y^\sigma)\}}+0.16m\le \left(\frac12+\frac13+0.16\right)m, 
		$$
		a contradiction.

		Thus, ${a(x^\sigma)=a(y^\sigma)=2}$. This implies that ${e_x=e_y=1}$; in the opposite case ${\Delta=1\bmod 8}$ by~\eqref{emoreover}, and one of $\Delta_x,\Delta_y$, being ${4\bmod 32}$,  cannot admit singular moduli with denominator~$2$, see  Proposition~\ref{pevenden}:\ref{i432none}.

		Since ${a(x^\sigma)=a(y^\sigma)=2}$ but ${x^\sigma\ne y^\sigma}$, the singular moduli~$x^\sigma$ and~$y^\sigma$  must be $4$-isogenous, see  Proposition~\ref{pvarisog}:\ref{isubdomisog}.  Hence so are~$x$ and~$y$. Corollary~\ref{cisog} now implies  that ${a_y=4}$, and Proposition~\ref{plinearel}  yields
		$$
		m'+\frac{n'}{4}-\frac{m'}{2}-\frac{n'}{2}=0.  
		$$
		Hence ${n'=2m'}$, again a contradiction.
		Thus,   ${n>0}$ implies that ${a(y^\sigma)=1}$. 
		
		Note that if~\eqref{emainoption} holds for some ${\sigma\in G}$ then it also holds with~$\sigma$ replaced by $\sigma^{-1}$. Hence ${n>0}$ implies that ${a(y^{\sigma^{-1}})=1}$ as well. Since there can be only one dominant singular modulus of given discriminant, we must have ${y^{\sigma}=y^{\sigma^{-1}}}$. Hence ${\sigma^2=1}$ by Lemma~\ref{lperm}. 
		
		\bigskip
		
		Now assume that ${n<0}$ but ${a(y) \ge 2}$. 
		The same argument as above shows that 
		${a(x^\sigma)=a(y)=2}$ and ${e_x=e_y=1}$. 
		
		The singular moduli $x^\sigma$ and~$y$ must be $4$-isogenous. Hence so are~$x$ and~$y^{\sigma^{-1}}$, which implies that ${a(y^{\sigma^{-1}})=4}$. Applying Proposition~\ref{pineq} with ${A=5}$ and ${\eps=0.16}$ to
		$$
		x^my^{-|n|}(x^{\sigma^{-1}})^{-m}(y^{\sigma^{-1}})^{|n|}=1,
		$$
		we obtain
		$$
		m'+\frac{|n'|}{4}\le 
		\frac{|n'|}{2}+ \frac{m'}{\min\{5,a(x^{\sigma^{-1}})\}}+0.16m'. 
		$$
		Since ${|n'|\le m'}$ and ${a(x^{\sigma^{-1}})\ge 2}$,  this is impossible.
		Thus, we proved that ${n<0}$ implies that ${a(y)=1}$. 
		
		\bigskip

		Finally,  ${(e_x,e_y)\ne (1,1)}$, because there cannot be two distinct dominant singular moduli of the same discriminant. Hence ${\Delta\equiv 1\bmod8}$ by~\eqref{emoreover}. 
		The proposition is proved.
	\end{proof}
	
	\subsubsection{Controlling the four denominators}
\label{sssfourdems}	
	Thus, we know that two of the singular moduli ${x,y,x^\sigma,y^\sigma}$ are dominant. Unfortunately, we have  no control over the denominators of the other two.
	
	We will show now that, with a suitably chosen Galois morphism~$\theta$, we can control the denominators of all four of 
	${x^\theta,y^\theta,x^{\sigma\theta},y^{\sigma\theta}}$.
	
Note that, so far, we assumed that ${\max\{|\Delta_x|, |\Delta_y|\}\ge 10^6}$ and used Proposition~\ref{pineq} with ${\eps=0.16}$. However, now we know that 
$$
\{\Delta_x,\Delta_y\}=\{\Delta, 4\Delta\} \quad \text{with}\quad \Delta\equiv 1\bmod 8,
$$	
which allows us (See Remark~\ref{repses}) to assume that  ${\max\{|\Delta_x|, |\Delta_y|\}\ge 10^8}$ and to use Proposition~\ref{pineq} with ${\eps=0.016}$.

	\begin{proposition}
		\label{ptwoeight}
		There exists ${\theta \in \Gal(L/K)}$ such that, when ${(e_x,e_y)=(1,2)}$, we have 
		\begin{equation}
		\label{etwoeight}
		a(x^\theta)=a(x^{\sigma\theta})=2, \qquad a(y^\theta)=a(y^{\sigma\theta})=8, 
		\end{equation}
		and when ${(e_x,e_y)=(2,1)}$, we have~\eqref{etwoeight} with~$x$ and~$y$ switched. 
	\end{proposition}
	
	To prove this proposition, we need  to bound ${|n'|}$ from below. 
	
	\begin{lemma}
		Assume that ${(e_x,e_y)=(2,1)}$. Then ${|n'|\ge 0.85m'}$. 
	\end{lemma}
	
	A similar estimate can be proved when ${(e_x,e_y)=(1,2)}$, but we do not need this.
	
	\begin{proof}
		Assume first that ${n<0}$. Then ${a(x)=a(y)=1}$. In particular,~$x$ and~$y$ are $2$-isogenous, and so are $x^\sigma,y^\sigma$. 
		Write 
		$$
		x^my^{-|n|}(x^\sigma)^{-m}(y^\sigma)^{|n|}=1.
		$$
		When ${a(x^\sigma)\ge 8}$ we use Proposition~\ref{pineq} with ${A=8}$ and ${\eps=0.016}$ to obtain 
		$$
		m'\le |n'|+\frac{m'}{8}+0.016m',
		$$
		which implies that ${|n'|\ge 0.85m'}$.
		
		When ${a(x^\sigma)\le 7}$, we must have ${a(x^\sigma)\in\{3,5,7\}}$ by  Proposition~\ref{pevenden}:\ref{i432none}, because  ${\Delta_x=4\Delta\equiv 4\bmod 32}$. Since~$x^\sigma$ and~$y^\sigma$ are $2$-isogenous, Corollary~\ref{cisog} implies that  ${a(y^\sigma)\in \{a(x^\sigma),a(x^\sigma)/4\}}$, and we must have ${a(y^\sigma)= a(x^\sigma)}$. Using Proposition~\ref{plinearel} with ${A=7}$ we obtain
		$$
		m'-|n'|-\frac {m'}{a(x^\sigma)}+\frac{|n'|}{a(x^\sigma)}=0, 
		$$
		which shows that ${|n'|=m'}$. This proves the lemma in the case ${n<0}$.

		Now assume that ${n>0}$. Then ${a(x)=a(y^\sigma)=1}$ and ${\sigma^2=1}$. In particular,~$x$ and~$y^\sigma$ are $2$-isogenous, and so are ${x^{\sigma^{-1}}=x^\sigma}$ and~$y$. Arguing as above, we obtain that either  ${a(x^\sigma)\ge 8}$ and ${n'\ge 0.85m'}$, or ${a(x^\sigma)\in\{3,5,7\}}$ and ${n'=m'}$. 
		The lemma is proved. 
	\end{proof}
	
	\begin{proof}[Proof of Proposition~\ref{ptwoeight}]

		Let us assume first that ${n<0}$ and  ${(e_x,e_y)=(1,2)}$.
		
		We have 
		\begin{equation}
		\label{edeltaxdeltay}
		\Delta_x=\Delta\equiv1\bmod8,  \qquad \Delta_y=4\Delta\equiv 4\bmod 32.  
		\end{equation}
		By Proposition~\ref{pevenden}:\ref{ionemodeight},  there exist two distinct morphisms ${\theta \in \Gal(L/K)}$ such that ${a(x^\theta)=2}$. Of the two, there can be  at most one with the property ${a(x^{\sigma\theta})=1}$. Hence we may find~$\theta$ satisfying 
		$$
		a(x^\theta)=2, \qquad a(x^{\sigma\theta})\ge 2. 
		$$
		Since ${n<0}$, we have ${a(y)=1}$ by Proposition~\ref{pinside}. Hence~$x$ and~$y$ are $2$-isogenous, and so are~$x^\theta$ and~$y^\theta$. It follows that ${a(y^\theta)\in \{2,8\}}$. But 
		${a(y^\theta)\ne 2}$ by  Proposition~\ref{pevenden}:\ref{i432none}. Hence 
		${a(y^\theta)=8}$.

		Proposition~\ref{pineq}, applied to 
		$$
		(x^\theta)^m(y^\theta)^{-|n|}(x^{\sigma\theta})^{-m}(y^{\sigma\theta})^{|n|}=1
		$$ 
		with ${A=9}$ and ${\eps=0.016}$, implies that 
		$$
		\frac{m'}2\le 
		\frac{|n'|}8+ \frac{m'}{\min\{9,a(x^{\sigma\theta})\}}+0.016m'. 
		$$
		If  ${a(x^{\sigma\theta})\ge 3}$ then this implies that 
		$$
		m'\left(\frac12-\frac13\right) \le \frac{|n'|}8+0.016m',
		$$
		which is impossible because ${m'\ge |n'|}$. Hence
		${a(x^{\sigma\theta})= 2}$, 
		and, as above, this implies that 
		${a(y^{\sigma\theta})= 8}$.

		\bigskip
		
		Now assume that ${n<0}$ and  ${(e_x,e_y)=(2,1)}$. We again have~\eqref{edeltaxdeltay}, but with~$x$ and~$y$ switched. 
		Arguing as before, we find  ${\theta \in \Gal(L/K)}$ such that 
		$$
		a(y^\theta)=2, \qquad a(y^{\sigma\theta})\ge 2, \qquad a(x^\theta)=8. 
		$$
		As before,  in the case ${a(y^{\sigma\theta})\ge 3}$ we apply  Proposition~\ref{pineq} to 
		$$
		(x^\theta)^{-m}(y^\theta)^{|n|}(x^{\sigma\theta})^{m}(y^{\sigma\theta})^{-|n|}=1
		$$ 
		and obtain 
		$$
		|n'|\left(\frac12-\frac13\right) \le \frac{m'}8+0.016m',
		$$
		which is impossible because ${|n'|\ge 0.85m'}$. Hence
		${a(y^{\sigma\theta})= 2}$, 
		and, as above, this implies that 
		${a(x^{\sigma\theta})= 8}$. 
		
		\bigskip

		Finally, let us assume that ${n>0}$. Then ${a(y^\sigma)=1}$ and ${\sigma^2=1}$. In particular,~$x$ and~$y^\sigma$ are $2$-isogenous, and so are~$x^\sigma$ and ${y^{\sigma^2}=y}$. Now, writing 
		$$
		x^m(y^\sigma)^{-n}(x^\sigma)^{-m}y^n=1, 
		$$
		we repeat the previous argument with $y,y^\sigma$ switched, and with~$n$ replaced by~$-n$. The proposition is proved. 
	\end{proof}

	\subsubsection{Completing the proof}
	\label{ssscomplete}
	
	Now we are ready to rule~\eqref{emainoption} out by deriving a contradiction. Let us summarize what we have. After renaming, we have distinct singular moduli $x_1,x_2$ of discriminant~$\Delta$  and $y_1,y_2$ of discriminant $4\Delta$ such that 
	$$
	a(x_1)=a(x_2)=2, \qquad a(y_1)=a(y_2)=8,  
	$$
	and
	\begin{equation}
	\label{eonetwo}
	x_1^{m_1}y_1^{n_1}x_2^{-m_1}y_2^{-n_1}=1, 
	\end{equation}
	where ${m_1,n_1}$ is a permutation of $m,n$. We want to show that this impossible.

	Proposition~\ref{pbasrelsigma} implies that we may assume 
	\begin{equation}
	\label{emaxmn}
	\max\{|m_1|,|n_1|\}\le 10^{10}|\Delta|^{1/2}. 
	\end{equation}
	Note also that 
	\begin{equation}
	\label{elowerdeltabefore}
	|\Delta|\ge 10^{7}
	\end{equation}
	by the assumption \eqref{elowerprimeleight}.

	Proposition~\ref{pevenden} implies that, after possible renumbering, we have
	\begin{align*}
	\tau(x_1)& =\frac{1+\sqrt\Delta}{4}, & \tau(x_2)&=\frac{-1+\sqrt\Delta}{4},\\
	\tau(y_1)&=\frac{b+\sqrt\Delta}{8},& \tau(y_2)&=\frac{-b+\sqrt\Delta}{8}, 
	\end{align*}
	where ${b\in \{\pm1,\pm3\}}$. 
	Denote ${t=e^{-\pi|\Delta|^{1/2}/4}}$ and ${\xi=e^{b\pi i/4}}$.  Then
	$$
	e^{2\pi i\tau(x_1)}= it^2, \quad e^{2\pi i\tau(x_2)}= -it^2, \quad e^{2\pi i\tau(y_1)}= \xi t, \quad 
	e^{2\pi i\tau(y_2)}= \bar\xi t.
	$$
	We deduce from~\eqref{eonetwo} that
	\begin{equation}
	\label{elogofrel}
	m_1\log(it^2x_1)-m_1\log(-it^2x_2)+n_1\log(\xi ty_1)-n_1\log(\bar\xi ty_2) \in \frac{1}{4}\pi i \Z.
	\end{equation}
Corollary~\ref{clogx}:\ref{icomplog} implies that
\begin{align*}
\log(it^2x_1)& =744it^2 + O_1(5\cdot10^5t^4), & \log(-it^2x_2)& =-744it^2 + O_1(5\cdot10^5t^4), \\
\log(\xi ty_1)& =744\xi t+ O_1(5\cdot10^5t^2), & \log(\bar\xi ty_2)& =744\bar\xi t+ O_1(5\cdot10^5t^2)
\end{align*}
Transforming the left-hand side of~\eqref{elogofrel} using  these expansions, we obtain 
	$$
	744(\xi-\bar\xi)tn_1+O_1\bigl(10^7t^2\max|m_1|,|n_1|\}\bigr) =\frac{1}{4}\pi i k  
	$$
	for some ${k\in \Z}$. 
	An easy estimate using~\eqref{emaxmn} and~\eqref{elowerdeltabefore} shows that the left-hand side does not exceed~$10^{-1000}$ in absolute value. Hence ${k=0}$, and we obtain, again using~\eqref{emaxmn} and~\eqref{elowerdeltabefore}, that 
	\begin{equation*}
	744|\xi-\bar\xi||n_1| \le 10^{17}|\Delta|^{1/2}e^{-\pi|\Delta|^{1/2}/4}< 10^{-900}. 
	\end{equation*}
	Hence ${n_1=0}$, a contradiction.

	This proves impossiblity of~\eqref{emainoption}, completing thereby the proof of Theorem~\ref{thprimel} in the case of equal fundamental discriminants.

	\subsection{Distinct fundamental discriminants}
	\label{ssdisfudis}
	In this subsection ${D_x\ne D_y}$. Arguing as in the beginning of Subsection~\ref{ssprepage}, we find a Galois morphism~$\sigma$ such that ${\Q(x,x^\sigma)=\Q(y,y^\sigma)}$. Lemma~\ref{lndeq}:\ref{inefudi}  implies that ${\Q(x)=\Q(y)}$. Corollary~\ref{cxy}  implies that ${h(\Delta_x) =h(\Delta_y) \le 16}$,  and our hypothesis 
	${\max\{|\Delta_x|, \Delta_y|\}\ge 10^6}$ contradicts Proposition~\ref{pwatki}.

	Theorem~\ref{thprimel} is proved.

	\section{Proof of Theorem~\ref{thpizthree}}
	\label{spizthree}
	In this section we prove Theorem~\ref{thpizthree}. 
	Throughout this section, unless the contrary is stated explicitly, 
	$x,y,z$ are distinct singular moduli satisfying 
	\begin{equation}
	\label{elower}
	\max\{|\Delta_x|,|\Delta_y|,|\Delta_z|\}\ge 10^{10},
	\end{equation}
	such that there exist non-zero integers ${m,n,r}$ with the property ${x^my^nz^r\in \Q^\times}$.

	We assume that $x,y,z$ have the same fundamental discriminant~$D$; if this is not the case, then the argument is much simpler, see Subsection~\ref{ssdis}.

	We denote by ${K=\Q(\sqrt D)}$ the common CM field of $x,y,z$, and we denote 
	$$
	L=K(x,y,z), \qquad G=\Gal(L/K). 
	$$
	We set
	$$
	f=\gcd(f_x,f_y,f_z), \quad e_x=\frac{f_x}f, \quad e_y=\frac{f_y}f, \quad e_z=\frac{f_z}f, \quad \Delta=Df^2.  
	$$ 
	Then ${\gcd(e_x,e_y,e_z)=1}$ and 
	\begin{align}
	\label{eexeyez}
	\Delta_x=e_x^2\Delta, \qquad \Delta_y=e_y^2\Delta, \qquad \Delta_z=e_z^2\Delta.
	\end{align}

	\subsection{The discriminants}
	The following property, showing that the ring class fields $K(x)$, $K(y)$ and $K(z)$ are closely related, 
is the basis for everything.
	
	\begin{proposition}
		\label{pexeyez}
		\begin{enumerate}
			\item
			\label{ikxykxzkyz}
			We have
			\begin{equation}
			\label{ekxykxzkyz}
			L=K(x,y)=K(x,z)=K(y,z).
			\end{equation}
			\item
			\label{ikxkykz}
			Each of the fields $K(x),K(y),K(z)$ is a subfield of~$L$ of degree at most~$2$. 
			
			\item
			\label{iexeyez}
			Up to permuting $x,y,z$ (and permuting correspondingly $m,n,r$) we have one of the cases from Table~\ref{taexeyez}. 
		\end{enumerate}
	\end{proposition}
	\begin{table}
		\caption{Data for Proposition~\ref{pexeyez}}
		\label{taexeyez}
		{\scriptsize
			$$
			\begin{array}{cccccclc}
			e_x&[L:K(x)]&e_y&[L:K(y)]&e_z&[L:K(z)]& \multicolumn{2}{c}{\text{remarks}}\\
%			\Delta\equiv&\text{remark}\\
			\hline
			1&
			1&
			1&
			1&
			1&
			1&\\
			1&
			1&
			1&
			1&
			2&
			1&
			\Delta\equiv1\bmod8\\
			1&
			1&
			2&
			1&
			2&
			1&
			\Delta\equiv1\bmod8\\
			1&
			2&
			2&
			1&
			2&
			1&
			\Delta\equiv0\bmod4,
			&n=r\\
			1&
			2&
			3&
			1&
			3&
			1&
			\Delta\equiv1\bmod3,
			&n=r\\
			2&
			2&
			3&
			1&
			3&
			1&
			\Delta\equiv1\bmod24,
			&n=r\\
			1&
			2&
			4&
			1&
			4&
			1&
			\Delta\equiv1\bmod8,&
			n=r\\
			1&
			2&
			6&
			1&
			6&
			1&
			\Delta\equiv1\bmod24,
			&n=r\\
			\end{array}
			$$}
	\end{table}

	\begin{proof}
		By the assumption, ${K(x^m)= K(y^nz^r)\subset K(y,z)}$. Lemma~\ref{lnorou} implies that ${K(x)= K(y^nz^r)}$, and, in particular, ${x\in  K(y,z)}$. Hence ${L=K(y,z)}$. By symmetry, we obtain~\eqref{ekxykxzkyz}.  This proves item~\ref{ikxykxzkyz}. 
		
		From~\eqref{elower} we may assume that, for instance, ${|\Delta_z|\ge 10^{10}}$. 
		Theorem~\ref{thprimel} implies that  the field ${K(x)= K(y^nz^r)}$   is a subfield of~$L$ of degree at most~$2$, and the same holds true for the fields $K(y)$.
		Unfortunately, we cannot make the same conclusion for   $K(z)$, because, \textit{a priori}, we cannot guarantee that ${\max\{|\Delta_x|,|\Delta_y|\}\ge 10^8}$, which is needed to apply Theorem~\ref{thprimel} in this case. So a rather lengthy extra argument is required to prove that ${[L:K(z)]\le 2}$. We split it into two cases. 
		
		Assume first that 
		\begin{equation}
		\label{egoodcases}
		\Delta\ne -3,-4 \quad \text{or} \quad \gcd(e_x,e_z) >1. 
		\end{equation}
		In this case, setting   ${\ell=\lcm(e_x,e_y,e_z)}$,   Proposition~\ref{pcompo} implies that ${L=K[\ell f]}$,  the ring class field of~$K$ of  conductor $\ell f$.  The class number formula~\eqref{eclnfr} implies that  
		\begin{equation}
		\label{ecnf}
		2\ge [L:K(x)] = \Psi(\ell/e_x, \Delta_x),  
		\end{equation}
		which results  in one of the following six cases:
		\begin{equation}
		\label{esixoptions}
		\begin{aligned}
		&\ell=e_x, && L=K(x);\\
		&\ell=2e_x, &&L=K(x), && \Delta_x\equiv 1\bmod 8;\\
		&\ell=2e_x, &&[L:K(x)]=2, && \Delta_x\equiv 0\bmod 4;\\
		&\ell=4e_x, &&[L:K(x)]=2, && \Delta_x\equiv 1\bmod 8;\\
		&\ell=3e_x, &&[L:K(x)]=2, && \Delta_x\equiv 1\bmod 3;\\
		&\ell=6e_x, &&[L:K(x)]=2, && \Delta_x\equiv 1\bmod 24. 
		\end{aligned}
		\end{equation}
		In particular, 
		${e_x\ge \ell/6\ge e_z/6}$, which implies that ${|\Delta_x|\ge |\Delta_z|/36\ge 10^8}$. Hence we may again apply Theorem~\ref{thprimel} to conclude that ${[L:K(z)]\le 2}$ in case~\eqref{egoodcases}.

		We are left with the case 
		\begin{equation}
		\label{eimpossi}
		\Delta \in \{-3,-4\} \quad\text{and}\quad \gcd(e_x,e_z)=1. 
		\end{equation} 
		We want to show  that it is impossible.  We claim that in this case ${\ph(e_z) \le 6}$.
		Indeed, if    ${x\in K}$ then
		$$
		2\ge [K(x,z):K(x)]=[K(z):K]= \frac{\Psi(e_z,\Delta)}{[\calO_K^\times:\calO_{K(z)}^\times]} \ge  \frac{\Psi(e_z,\Delta)}{3}, 
		$$
		which proves that ${\ph(e_z)\le \Psi(e_z, \Delta)\le 6}$. And if ${x\notin K}$ then 
		$$
		\Psi(e_z, e_x^2\Delta) = \bigl [K[e_ze_x]:K[e_x]\bigr]=  \bigl [K[e_ze_x]:K(x,z)\bigr]\cdot [K(x,z) :K(x)]. 
		$$
		We have ${\bigl [K[e_ze_x]:K(x,z)\bigr]\le 3}$  by Proposition~\ref{pcompo}, and 
		$$
		[K(x,z) :K(x)]=[L:K(x)]\le 2,
		$$
		as we have seen above. It follows that ${\ph(e_z) \le \Psi(e_z, e_x^2\Delta)\le 6}$.

		From 
		${\ph(e_z) \le 6}$
		we deduce ${e_z\le 18}$. Hence ${|\Delta_z|\le 4\cdot18^2<10^{10}}$, a contradiction. This shows impossibility of~\eqref{eimpossi}, completing thereby the proof of item~\ref{ikxkykz}. 
		
		We are left with part~\ref{iexeyez} of the proposition. 
		As we have just seen, we have one of the six cases~\eqref{esixoptions}, and similarly with~$x$ replaced by~$y$ or by~$z$. 
Since ${e_x,e_y,e_z}$ are coprime, every prime number~$p$ does not divide one of them. If, say, ${p\nmid e_x}$, then $$
\nu_p(\ell) =\nu_p(\ell/e_x) \le 
\begin{cases}
2, &p=2, \\
1, & p=3, \\
0, &p\ge 5. 
\end{cases}
$$
This proves that ${\ell\mid 12}$. 		
Moreover, 
		\begin{align}
		\label{etwo}
		\text{if ${2\mid \ell}$\ } &\text{then ${\Delta\equiv 1\bmod 8}$ or ${\Delta\equiv 0\bmod 4}$}; \\
		\label{efour}
		\text{if ${4\mid \ell}$\ } &\text{then ${\Delta\equiv 1\bmod 8}$}; \\
		\label{ethree}
		\text{if ${3\mid \ell}$\ } &\text{then ${\Delta\equiv 1\bmod 3}$}. 
		\end{align}
		Indeed, assume that  ${2\mid \ell}$ but ${\Delta\equiv 5\bmod 8}$. Since ${e_x,e_y,e_z}$ are coprime, one of them, say,~$e_x$, is not divisible by~$2$. Then ${(\Delta_x/2)=-1}$, and ${\Psi(\ell/e_x, \Delta_x)}$  must be divisible by~$3$, which contradicts~\eqref{ecnf}.  This proves~\eqref{etwo}.
		
		In a similar fashion, one shows that 
		${\Psi(\ell/e_x, \Delta_x)}$ is divisible by~$4$ in each of the cases
		\begin{align*}
		&4\mid \ell, \qquad  2\nmid e_x, \qquad \Delta \equiv 0\bmod 4,\\
		&3\mid \ell, \qquad  3\nmid e_x, \qquad \Delta \equiv 2\bmod 3,
		\end{align*}
		and  it is divisible by~$3$ in the case
		$$
		3\mid \ell, \qquad  3\nmid e_x, \qquad \Delta \equiv 0\bmod 3. 
		$$
		This proves~\eqref{efour} and~\eqref{ethree}. 
		
		Finally, it also follows from Theorem~\ref{thprimel} that, when, say, ${K(x)\ne L}$, we must have 
		${e_y=e_z}$ and ${n=r}$. 
		
		A little \textsf{PARI} script  (or verification by hand) shows that, up to permuting $x,y,z$, all possible cases are listed in Table~\ref{taexeyez}.
	\end{proof}

	Recall that we denote ${f=\gcd(f_x,f_y,f_z)}$. Denote ${L_0=K[f]}$. Recall also that ${G=\Gal(L/K)}$.  Then we have the following consequence.
	
	\begin{corollary}
		\label{cocases}
		\begin{enumerate}
			\item
			\label{ioneortwo}
			We have either 
			$$
			L=L_0=K(x)=K(y)=K(z),
			$$
			or ${[L:L_0]=2}$, in which case exactly one of the fields $K(x),K(y),K(z)$ is~$L_0$ and the other two are~$L$.
			
			\item
			\label{ihundredone}
			We have ${[L_0:K]\ge 101}$.
			
			\item
			\label{ibigdensthree}
			There exists ${\sigma \in G}$ such that 
			${a(x^\sigma), a(y^\sigma),a(z^\sigma)\ge 13}$, 
			
			\item
			\label{ibigdenstwo}
			There exists ${\sigma \in G}$ such that 
			${a(x^\sigma), a(y^\sigma)\ge 18}$, 
			and the same statement holds true for $x,z$ and for $y,z$.

			\item
			\label{ibigdenone}
			There exists ${\sigma \in G}$ such that 
			${a(x^\sigma)\ge 30}$, 
			and the same statement holds true for~$y$ and for~$z$. 
			
		\end{enumerate}
		
	\end{corollary}
	
	\begin{proof}
		Item~\ref{ioneortwo} is proved just by exploring Table~\ref{taexeyez}. To prove item~\ref{ihundredone}, note that, since ${\max\{e_x,e_y,e_z\}\le 6}$, we have 
		\begin{equation}
		\label{elowerdelta}
		|\Delta|\ge \max\{|\Delta_x|,|\Delta_y|,|\Delta_z|\}/36\ge  10^8
		\end{equation}
		by~\eqref{elower}. Hence 
		${[L_0:K] =h(\Delta) \ge 101}$
		by Proposition~\ref{pwatki}. 
		
		In proving item~\ref{ibigdensthree} we must distinguish the case ${L=L_0}$ and ${[L:L_0]=2}$. In the former case 
		${L=K(x)=K(y)=K(z)}$ and ${|G|=[L:K]\ge 101}$. Proposition~\ref{pcountden} implies that there exist at most~$32$ elements ${\sigma\in G}$ such that ${a(x^\sigma)< 13}$, and the same for~$y$ and~$z$. Since ${|G|\ge 101>96}$, we can find ${\sigma \in G}$ as wanted. 
		
		If ${[L:L_0]=2}$ then, say, 
		$$
		K(x)=L_0, \qquad K(y)= K(z)=L,
		$$
		and
		${|G|=[L:K]\ge 202}$. Again using Proposition~\ref{pcountden}, there exist at most~$32$ elements ${\sigma\in G}$ such that ${a(y^\sigma)<13}$, the same for~$z$, and at most~$64$ elements ${\sigma\in G}$ such that ${a(x^\sigma)<13}$. Since ${|G|\ge 202>128}$, we again can find a~$\sigma$ as wanted. This proves item~\ref{ibigdensthree}.
		
		Item~\ref{ibigdenstwo} is proved similarly. In the case ${L=K(x)=K(y)}$ there exist at most~$48$ elements ${\sigma\in G}$ such that ${a(x^\sigma)< 18}$, and the same for~$y$. Since ${|G|\ge 101>96}$, we are done. In the case when one of $K(x),K(y)$ is~$L$, the other is~$L_0$, and 
		${[L:L_0]=2}$, we have ${48+96=144}$ unsuitable ${\sigma\in G}$; since ${|G|\ge 202}$, we are done again. 
		
		Item~\ref{ibigdenone} is similar as well: there exist at most~$99$ unsuitable~$\sigma$ when ${L=K(x)}$, and at most $198$ when ${[L:K(x)]=2}$; in both cases we conclude as before. 
	\end{proof}

	In the sequel we denote 
	$$
	m'=me_x, \qquad n'=ne_y, \qquad r'=re_z.
	$$
	We may and will assume that ${m>0}$ and that 
	\begin{equation}
	\label{empgenprp}
	m'\ge \max\{|n'|,|r'|\}, \qquad a_x=1. 
	\end{equation}
	In the course of the argument we will study  multiplicative relations
	\begin{equation}
	\label{emainmultrelationthree}
	x^my^nz^r(x^\sigma)^{-m}(y^\sigma)^{-n}(z^\sigma)^{-r}=1,
	\end{equation}
	with various choices of ${\sigma \in G}$, 
	using Propositions~\ref{plinearel} and~\ref{pineq}. 
	In our usage of Propositions~\ref{plinearel} and~\ref{pineq}  the  parameters therein will satisfy the following restrictions:
	\begin{equation}
	\label{eparametersthree}
	\begin{aligned}
	&k= 6, \quad X=\max\{|\Delta_x|,\Delta_y|,|\Delta_z|\}\ge 10^{10}, \quad Y=|\Delta|\ge \frac1{36}X, \\
	& A \le 162 \quad \text{for Proposition~\ref{plinearel}}, \\
	&A\le 30, \quad \eps=0.01 \quad \text{for Proposition~\ref{pineq}}. 
	\end{aligned} 
	\end{equation}
	It is easy to verify that for any choice of parameters satisfying~\eqref{eparametersthree}, conditions~\eqref{erootofy} and~\eqref{eassumpdelta} are met, so using the propositions is justified.

	\subsection{The denominators}
	We already know that~$x$ is dominant, see~\eqref{empgenprp}.  Our principal observation is that either one of $y,z$ is dominant as well, or they both are subdominant. More precisely, we have the following.

	\begin{proposition}
		\label{pdomorsubd}
Up to interchanging~$y$ and~$z$,  one of the following conditions is satisfied:
		\begin{align}
		\label{edomopt}
		&a_y=1\hphantom{=a_z=2}  \text{and}\quad n<0;\\
		\label{esubdopt}
		&a_y=a_z=2\hphantom{=1}  \text{and}\quad n,r<0. 
		\end{align}
\end{proposition}
	
	\begin{proof}
		With~$y$ and~$z$ possibly switched, we may assume that we are in one of the following cases:
		\begin{align}
		\label{eplusplus}
		&n,r>0;\\
		\label{eplusminus}
		&n<0 \quad\text{and}\quad r>0;\\
		\label{eminusminus}
		&n,r<0. 
		\end{align}
		We consider them separately. 
		
		\bigskip
		
		Assume~\eqref{eplusplus}. Let~$\sigma$ be like in  Corollary~\ref{cocases}:\ref{ibigdensthree}.  
		Applying Proposition~\ref{pineq} to
		$$
		x^my^nz^r(x^\sigma)^{-m}(y^\sigma)^{-n}(z^\sigma)^{-r}=1
		$$
		with ${A=13}$ and ${\eps=0.01}$, 
		and using that ${\max\{|m'|,|n'|,|r'|\}=m'}$ by~\eqref{empgenprp}, we obtain 
		\begin{equation*}
		m' \le \frac{m'}{13}+\frac{n'}{13}+\frac{r'}{13}+0.01m'\le m'\left(\frac3{13}+0.01\right), 
		\end{equation*}
		a contradiction. This shows that~\eqref{eplusplus} is impossible.
		
		\bigskip
		
		Now assume~\eqref{eplusminus}. We want to show that ${a_y=1}$ in this case. Thus, assume that ${a_y\ge 2}$. Using  Corollary~\ref{cocases}:\ref{ibigdenstwo}, we find ${\sigma \in G}$ such that ${a(x^\sigma),a(z^\sigma)\ge 18}$. 
		Applying Proposition~\ref{pineq} to
		$$
		x^my^{-|n|}z^r(x^\sigma)^{-m}(y^\sigma)^{|n|}(z^\sigma)^{-r}=1
		$$
		with ${A=18}$ and ${\eps=0.01}$, we obtain 
		\begin{equation*}
		m' \le \frac{|n'|}{\min\{18,a(y)\}}+\frac{m'}{18}+\frac{r'}{18}+0.01m'\le m'\left(\frac12+\frac2{18}+0.01\right), 
		\end{equation*}
		a contradiction. This shows that in the case~\eqref{eplusminus} we must have~\eqref{edomopt}. 
		
		\bigskip
		
		Finally, let us assume~\eqref{eminusminus}. If one of $a_y,a_z$ is~$1$ then we have~\eqref{edomopt}, possibly after switching. If ${a_y=a_z=2}$ then we have~\eqref{esubdopt}. Now let us assume that none of these is the case; that is,  both $a_y,a_z$ are $\ge2$ and one of them is $\ge3$.  
		Again using Corollary~\ref{cocases}, we may find ${\sigma \in G}$ such that ${a(x^\sigma)\ge 30}$. Applying Proposition~\ref{pineq}, we obtain, in the same fashion as in the previous cases, the inequality
		$$
		m'\le m'\left(\frac12+\frac13+\frac1{30}+0.01\right), 
		$$ 
		a contradiction. 
		The proposition is proved. 
	\end{proof}

	We study cases~\eqref{edomopt} and~\eqref{esubdopt} in Subsections~\ref{ssdomopt} and~\ref{sssubdopt}, respectively.

	\subsection{The dominant case}
	\label{ssdomopt}
	
	In this subsection we assume~\eqref{edomopt}. Thus, we have the following:
	$$
	m>0, \quad n<0, \quad m'\ge\max\{ |n'|,|r'|\}, \qquad a_x=a_y=1. 
	$$
	Since both~$x$ and~$y$ are dominant, we must have ${e_x\ne e_y}$. Exploring Table~\ref{taexeyez}, we find ourselves in one of the following cases:
	\begin{align}
	\label{ecaseonetwoeight}
	&\{e_x,e_y\}=\{1,2\},&& e_z\in\{1,2\}, && \Delta\equiv 1\bmod 8,\\
	\label{ecaseonetwotwofour}
	&\{e_x,e_y\}=\{1,2\}, && e_z=2, && \Delta\equiv 0\bmod 4,\\
	\label{ecaseonethreethreethree}
	&\{e_x,e_y\}=\{1,3\},&& e_z=3, && \Delta\equiv 1\bmod 3,\\
	\label{ecasetwothreethreetwentyfour}
	&\{e_x,e_y\}=\{2,3\},&& e_z=3, && \Delta\equiv 1\bmod 24,\\
	\label{ecaseonefourfoureight}
	&\{e_x,e_y\}=\{1,4\},&& e_z=4, && \Delta\equiv 1\bmod 8,\\
	\label{ecaseonesixsixtwentyfour}
	&\{e_x,e_y\}=\{1,6\},&& e_z=6, && \Delta\equiv 1\bmod 24. 
	\end{align}

	\begin{remark}
		It is absolutely crucial that, in each of the cases above, a non-trivial congruence condition is imposed on~$\Delta$. This would allow us to use Propositions~\ref{poddpden} and~\ref{pevenden} to find Galois morphisms~$\sigma$ with well-controlled denominators of $x^\sigma,y^\sigma,z^\sigma$, which is needed for the strategy described in Subsection~\ref{sssstra}  to work.  
	\end{remark} 
	
	Here are some more specific observations.
	
	\begin{enumerate}
		\item
		We have either ${e_z=e_x}$ or ${e_z=e_y}$, which implies, in particular, that 
		\begin{equation}
		\label{eazneone}
		a_z\ne 1. 
		\end{equation}

		\item
		In case~\eqref{ecaseonetwoeight}  we have ${K(x)=K(y)=K(z)=L}$.
		
		\item
		In cases \eqref{ecaseonetwotwofour}--\eqref{ecaseonesixsixtwentyfour} we have ${K(z)=L}$, and  one of the fields $K(x)$ or~$K(y)$ is~$L$ as well, while the other  is a degree~$2$ subfield of~$L$. More precisely:
		\begin{itemize}
			\item
			if ${e_x<e_y=e_z}$ then ${K(y)=L}$ and ${[L:K(x)]=2}$;
			\item
			if ${e_y<e_x=e_z}$ then ${K(x)=L}$ and ${[L:K(y)]=2}$. 
		\end{itemize}

		\item
		Theorem~\ref{thprimel} implies  that in 
		cases \eqref{ecaseonetwotwofour}--\eqref{ecaseonesixsixtwentyfour} we have
		${n=r}$ when ${e_x<e_y}$, and ${m=r}$ when ${e_x>e_y}$.
		
	\end{enumerate}

	\subsubsection{The strategy}
	\label{sssstra}
	
	In each of  cases \eqref{ecaseonetwoeight}--\eqref{ecaseonesixsixtwentyfour} we apply the following strategy.
	
	\begin{itemize}
		\item
		Find possible values for~$a_z$.
		
		\item
		Using Proposition~\ref{poddpden} or~\ref{pevenden},  find several ${\sigma\in G}$ such that we can control the denominators 
		\begin{equation}
		\label{ethreedens}
		a(x^\sigma), \qquad a(y^\sigma), \qquad a(z^\sigma). 
		\end{equation}

		\item
		For every such~$\sigma$, and every possible choice of~$a_z$ and of denominators~\eqref{ethreedens}, Proposition~\ref{plinearel} implies the linear equation
		$$
		m'+n' + \frac{r'}{a_z} = \frac{m'}{a(x^{\sigma})}+ \frac{n'}{a(y^{\sigma})}+ \frac{r'}{a(z^{\sigma})}.
		$$  
		With sufficiently many choices of~$\sigma$, we may hope to have sufficiently many equations to conclude that ${m'=n'=r'=0}$, a contradiction. 
	\end{itemize}
	
	Practical implementation of this strategy differs from case to case. For instance, in cases~\eqref{ecaseonetwotwofour}--\eqref{ecaseonesixsixtwentyfour} we have ${m'=r'}$ or ${n'=r'}$, so wee need only two independent equations to succeed, while in case~\eqref{ecaseonetwoeight}  three independent equations are needed.
	
	Case~\eqref{ecaseonetwotwofour} is somewhat special, because we get only one equation. To complete the proof in that case, we need  to use an argument similar to that of Subsection~\ref{ssscomplete}. 
	
	Below  details for all the cases follow.

	\subsubsection{Cases~\eqref{ecaseonethreethreethree}--\eqref{ecaseonesixsixtwentyfour}}
	\label{sssfourcases}
	
	In these cases ${K(z)=L}$, and one of the fields~$K(x)$, $K(y)$ is also~$L$ while the other is a degree~$2$ subfield of~$L$ . In this subsection we make no use of the assumption ${m'\ge |n'|}$. Hence we may assume that ${e_x<e_y=e_z}$, in which case we have 
	\begin{equation}
	\label{eyzxdegs}
	K(y)=K(z)=L, \qquad [L:K(x)]=2.
	\end{equation}
	Moreover, Theorem~\ref{thprimel} implies that in this case ${n=r}$,  and that $y,z$ are conjugate over $K(x)$:
	\begin{equation}
	\label{exthetay}
	y^\theta=z, \qquad z^\theta=y, 
	\end{equation}
	where~$\theta$ is the non-trivial element of $\Gal(L/K(x))$. 
	
	Let us specify the general strategy described in Subsection~\ref{sssstra} for the cases~\eqref{ecaseonethreethreethree}--\eqref{ecaseonesixsixtwentyfour}.

	\begin{enumerate}
		\item
		\label{iaz}
		Proposition~\ref{pvarisog} implies that~$x$ and~$y$ are $\ell$-isogenous, where ${\ell=e_xe_y}$. 
		Hence ${x=x^\theta}$ and ${z=y^\theta}$ are $\ell$-isogenous as well. Using Corollary~\ref{cisog}, we may now shortlist possible values of the denominator $a_z$. Precisely, 
		$$
		a_z \in \left(\frac{e_z}{e_x}\calQ(\ell)\right) \cap\Z_{\ge 2},  
		$$
		where we use notation ${\lambda S=\{\lambda s: s\in S\}}$.
		For instance, in  case~\eqref{ecasetwothreethreetwentyfour} we have ${\ell=6}$, and
		$$
		a_z\in \left(\frac32\left\{\frac16,\frac23,\frac32,6\right\}\right)\cap\Z_{\ge 2} = \{9\}. 
		$$
		
		\item
		Propositions~\ref{poddpden}:\ref{ionemodthree} and~\ref{pevenden}:\ref{ionemodeight} imply existence of morphisms~$\sigma_1$ and~$\sigma_2$ such that the three denominators 
		$a(x)$ (which is~$1$), $a(x^{\sigma_1})$ and $a(x^{\sigma_2})$ 
		are distinct. Precisely, 
		\begin{itemize}
			\item
			if ${\Delta_x\equiv1\bmod 3}$ then~$3$ and~$9$ are denominators for~$\Delta_x$; 
			\item
			if ${\Delta_x\equiv1\bmod 8}$ then~$2$ and~$4$ are denominators for~$\Delta_x$. 
		\end{itemize}
		
		For instance, in  case~\eqref{ecasetwothreethreetwentyfour} we may find~$\sigma_1$ and~$\sigma_2$ to have 
		$$
		a(x^{\sigma_1})=3, \qquad a(x^{\sigma_2})=9.
		$$

		\item
		Using again Corollary~\ref{cisog}, we may now shortlist the denominators $a(y^{\sigma_i})$ and $a(z^{\sigma_i})$. Precisely, 
		$$
		a(y^{\sigma_i}), a(z^{\sigma_i}) \in \left(a(x^{\sigma_i})\frac{e_z}{e_x}\calQ(\ell)\right) \cap\Z_{\ge 1}.   
		$$
		For instance, in  case~\eqref{ecasetwothreethreetwentyfour} we have 
		\begin{align*}
		a(y^{\sigma_1}), a(z^{\sigma_1}) &\in \left(3\cdot\frac32\left\{\frac16,\frac23,\frac32,6\right\}\right)\cap\Z_{\ge 1} = \{3,27\}, \\ 
		a(y^{\sigma_2}), a(z^{\sigma_2}) &\in \left(9\cdot\frac32\left\{\frac16,\frac23,\frac32,6\right\}\right)\cap\Z_{\ge 1} = \{9,81\}.
		\end{align*}

		\item
		Now, Proposition~\ref{plinearel} implies the system of linear equations 
		\begin{equation}
		\label{eequations4}
		m'+\left(1 + \frac{1}{a(z)}\right)n' = \frac{m'}{a(x^{\sigma_i})}+ \left(\frac{1}{a(y^{\sigma_i})}+ \frac{1}{a(z^{\sigma_i})}\right)n'\qquad (i=1,2).
		\end{equation}
		(recall that ${n'=r'}$). 
		Solving the system, we find that ${m'=n'=0}$ in every instance, a contradiction. This shows impossibility of cases~\eqref{ecaseonethreethreethree}--\eqref{ecaseonesixsixtwentyfour}.

		For instance, in  case~\eqref{ecasetwothreethreetwentyfour}, 
		equations~\eqref{eequations4} become 
		\begin{align*}
		m'+\left(1+\frac19\right)n'&= \frac{m'}{3}+ \lambda n', \qquad \lambda\in \left\{\frac23,\frac13+\frac{1}{27}, \frac{2}{27}\right\}\\
		m'+\left(1+\frac19\right)n'&= \frac{m'}{9}+\mu n', \qquad \mu\in \left\{\frac29,\frac19+\frac{1}{81}, \frac{2}{81}\right\}, 
		\end{align*}
		so nine systems in total, each of them having ${m'=n'=0}$ as the only solution.

	\end{enumerate}

	The numerical data obtained following these steps can be found in  Table~\ref{tafourcases}. Note that we have 390 linear systems to solve: nine systems in  cases \eqref{ecaseonethreethreethree}, \eqref{ecasetwothreethreetwentyfour}, 72 systems in case~\eqref{ecaseonefourfoureight}, and 300 systems in case \eqref{ecaseonesixsixtwentyfour}. Doing this by hand is impractical, and we used a \textsf{PARI} script for composing Table~\ref{tafourcases} and solving the systems.

	\begin{table}
		\caption{Cases~\eqref{ecaseonethreethreethree}--\eqref{ecaseonesixsixtwentyfour}  and case~\eqref{ecaseonetwotwofour} with ${\Delta\equiv 4\bmod 32}$}
		\label{tafourcases}
		{\tiny
			$$
			\begin{array}{l|lc@{\ }ccc@{}c@{}c@{}c@{}c@{}c}
			\text{case}&\Delta\equiv&e_x&e_y,&\ell& a_z&a(x^{\sigma_1})&a(y^{\sigma_1}), &a(x^{\sigma_2})&a(y^{\sigma_2}),&\text{total} \\
			&&&e_z&& &&a(z^{\sigma_1}) &&a(z^{\sigma_2})&\text{systems} \\
			\hline
			\text{\eqref{ecaseonethreethreethree}}&1\bmod 3&1&3&3&9&3&\in\{3,27\}&9&\in\{9,81\}&9\\
			\text{\eqref{ecasetwothreethreetwentyfour}}&1\bmod 24&2&3&6&9&3&\in\{3,27\}&9&\in\{9,81\}&9\\
			\text{\eqref{ecaseonefourfoureight}}&1\bmod 8&1&4&4&\in\{4,16\}&2&\in\{2,8,32\}&4&\in\{4,16,64\}&72\\
			\text{\eqref{ecaseonesixsixtwentyfour}}&1\bmod 24&1&6&6&\in\{4,9,36\}&2&\in\{2,8,18,72\}&3&\in\{3,12,27,108\}&300\\
			\text{\eqref{ecaseonetwotwofour}}&4\bmod 32&1&2&2&4&8&\in \{8,32\}&16&\in \{16,64\}&9
			\end{array}
			$$
		}
	\end{table}
	
	\begin{remark}
		Using  Propositions~\ref{poddpden} and~\ref{pevenden}, we can further refine the lists of possible denominators for~$z$,~$y^{\sigma_i}$ and~$z^{\sigma_i}$. For instance, 
		if the discriminant ${\Delta_y=\Delta_z}$ is ${\equiv0\bmod 9}$ then it cannot have denominators divisible by~$3$ but not by~$9$. Thus, in  case~\eqref{ecasetwothreethreetwentyfour}, number~$3$ cannot be the denominator of~$y^{\sigma_1}$ or of~$z^{\sigma_1}$, and so we must have ${a(y^{\sigma_1})= a(z^{\sigma_1})=27}$. Arguments of this kind, used systematically, allow one to decimate the number of systems to solve.
		
		However,  the computational time for solving our systems being insignificant, we prefer to disregard this observation. 
	\end{remark}

	\subsubsection{Case~\eqref{ecaseonetwotwofour}}
	\label{sssonetwotwofour}
	This case is similar to cases~\eqref{ecaseonethreethreethree}--\eqref{ecaseonesixsixtwentyfour}, but somewhat special. Let us reproduce our data for the reader's convenience:
	$$
	\{e_x,e_y\}=\{1,2\}, \quad e_z=2, \quad  \Delta\equiv 0\bmod 4, \quad  a_x=a_y=1. 
	$$ 
	We may again assume that ${e_x<e_y}$, which means now that ${e_x=1}$ and ${e_y=2}$, and we again have~\eqref{eyzxdegs},~\eqref{exthetay}. Furthermore, step~\ref{iaz} of the strategy described in Subsection~\ref{sssfourcases} works here as well: we prove that each of~$y,z$ is $2$-isogenous to~$x$, which, in particular,  allows us to determine that
	${a_z=4}$. For later use, let us note that 
	\begin{equation}
	\label{etauxyz}
	\tau_x=\frac{\sqrt\Delta}{2}, \qquad \tau_y=\sqrt\Delta, \qquad \tau_z= \frac{b'+\sqrt\Delta}{4}  ,  \quad 
	\end{equation}
	where 
	\begin{equation}
	\label{ebprimebis}
	b'=
	\begin{cases}
	0, &\text{if $\Delta\equiv 4\bmod 8$}, \\
	2, &\text{if $\Delta\equiv 0\bmod 8$}. 
	\end{cases}
	\end{equation}

	The rest of the argument splits into two subcases. If ${\Delta\equiv 4 \bmod 32}$ then we may proceed as in Subsection~\ref{sssfourcases}. Proposition~\ref{pevenden} implies that there exist ${\sigma_1,\sigma_2\in G}$ such that 
	${a(x^{\sigma_1})=8}$ and ${a(x^{\sigma_2})=16}$. 
	As before, we can now determine possible denominators of $y^{\sigma_i}, z^{\sigma_i}$ (see the bottom line of Table~\ref{tafourcases}) and solve the resulting systems~\eqref{eequations4}, concluding that ${m=n=0}$.

	Now assume that ${\Delta\not\equiv 4 \bmod 32}$. In this case~$2$ or~$4$ is a denominator for~$\Delta$, see  Proposition~\ref{pevenden}:\ref{ieveryn}.  Since ${\Delta_x=\Delta}$, there exists ${\sigma \in G}$ such that ${a(x^\sigma) \in \{2,4\}}$. We claim that 
	\begin{equation}
	\label{ealldiff}
	y^\sigma,z^\sigma \notin\{ y,z\}. 
	\end{equation}
	Indeed, if, say, ${y^\sigma=y}$ then ${\sigma = \id}$ because ${L=K(y)}$; but ${x^\sigma\ne x}$, a contradiction. For the same reason, ${z^\sigma\ne z}$. Now assume that ${y^\sigma=z}$. Theorem~\ref{thprimel} implies that~$y$ and~$z$ are conjugate over $K(x)$. Hence there exists ${\theta \in G}$ such that 
	$$
	x^\theta=x, \qquad y^\theta=z, \qquad z^\theta =y. 
	$$
	Then ${z^{\theta\sigma}=z}$, and, as before, ${\theta\sigma=\id}$, which is again a contradiction because ${x^{\theta\sigma}=x^\sigma\ne x}$. Similarly one shows that ${z^\sigma\ne y}$. This proves~\eqref{ealldiff}.

	The cases ${a(x^\sigma)=2}$ and ${a(x^\sigma)=4}$ are very similar, but each one has some peculiarities, so we consider them separately. 
	
	\bigskip
	
	Assume that ${a(x^\sigma)=2}$. Then ${a(y^\sigma)=a(z^\sigma)=8}$. Proposition~\ref{plinearel} gives
	$$
	m'+n'\left(1+\frac{1}{4}\right)=\frac{m'}{2}+n'\left(\frac18+\frac18\right),
	$$
	which is just ${m'=-2n'}$. Hence ${m=-4n}$. It follows that 
	${(x/x^\sigma)^4(y^\sigma/y)(z^\sigma/z)}$ is a root of unity. Since the roots of unity in~$L$ are of order dividing~$24$ (see Corollary~\ref{crofoneinrcf}), we obtain 
	\begin{equation}
	\label{e96}
	\bigl(x^{4}(x^\sigma)^{-4}y^{-1}z^{-1}y^\sigma z^\sigma\bigr)^{24}=1. 
	\end{equation}
	Now we are going to argue as in Subsection~\ref{ssscomplete}. This means:
	\begin{itemize}
		\item
		We  give explicit expressions for the $\tau$- and $q$-parameters of all the six singular moduli occurring in~\eqref{e96}.   Note that the $\tau$-parameters for $x,y,z$ are already given in~\eqref{etauxyz}, so we need to determine them only for $x^\sigma,y^\sigma,z^\sigma$. 
		
		\item
		Taking the logarithm of~\eqref{e96}, we deduce that a certain linear combination of logarithms is a multiple of $\pi i/12$. 
		
		\item
		Using the $q$-expansion from Corollary~\ref{clogx}, we obtain a contradiction. 
		
	\end{itemize}
	Note however that in Subsection~\ref{ssscomplete} the first order expansion~\eqref{efirst} was sufficient, while now we would need the second order expansion~\eqref{esecond}.

	Since $y,z,x^\sigma$ are distinct and $2$-isogenous to~$x$, we must have, in addition to~\eqref{etauxyz},~\eqref{ebprimebis}, 
	\begin{equation}
	\label{ebsecond}
	\tau(x^\sigma)= \frac{b_2+\sqrt\Delta}{4}, \quad \text{where} \quad 
	b_2=
	\begin{cases}
	0, &\text{if $b'=2$}, \\
	2, &\text{if $b'=0$}.
	\end{cases}
	\end{equation}
	Furthermore, since $x,y^\sigma, z^\sigma$ are distinct and $2$-isogenous to~$x^\sigma$, we must have 
	$$
	\{\tau(y^\sigma), \tau(z^\sigma)\} =\left\{\frac{b_2+\sqrt\Delta}{8},\frac{b_2'+\sqrt\Delta}{8}\right\},  \qquad 
	b_2' \in \{b_2+4,b_2-4\}. 
	$$
	Denote ${t=e^{-\pi|\Delta|^{1/2}/4}}$ and ${\xi=e^{\pi ib_2/4}}$.  Note that ${\xi\in \{1,i\}}$, and that 
	$$
	e^{\pi ib'/2}= - \xi^2, \qquad e^{\pi ib_2'/4}=-\xi. 
	$$
	We obtain
	\begin{align*}
	&e^{2\pi i\tau_x}=t^4, \qquad e^{2\pi i \tau_y} =t^8, \qquad e^{2\pi i \tau_z} = -\xi^2t^2, \qquad 
	e^{2\pi i \tau(x^\sigma)} = \xi^2t^2,\\
	& \bigl\{e^{2\pi i \tau(y^\sigma)}, e^{2\pi i \tau(z^\sigma)}\bigr\} =\{\xi t,-\xi t\}. 
	\end{align*}
	Taking the logarithm of~\eqref{e96}, we obtain
	$$
	4 \log(t^4x)-4\log (\xi^2t^2x^\sigma) -\log(t^8y)-\log(-\xi^2t^2z)+\log(-\xi t \cdot \xi t \cdot y^\sigma\cdot z^\sigma) \in \frac{\pi i}{12}\Z. 
	$$
	The $q$-expansion~\eqref{esecond} from Corollary~\ref{clogx} implies that for some ${k\in \Z}$ we have 
	$$
	\frac{\pi i}{12}k=-162000\xi^2t^2+O_1(10^{10}t^3). 
	$$
	This easily leads to contradiction, exactly like in Subsection~\ref{ssscomplete}. 

	\bigskip
	
	Now assume  that ${a(x^\sigma)=4}$. Since~$x^\sigma$ is $2$-isogenous to~$y^\sigma$ and to~$z^\sigma$, Corollary~\ref{cisog} and Proposition~\ref{pevenden} imply that ${a(y^\sigma),a(z^\sigma) \in \{4,16\}}$. Note however that 
	$$
	\Delta_y=\Delta_z=4\Delta\equiv 0\bmod 16,
	$$ 
	and   Proposition~\ref{pevenden}:\ref{i16none} implies that there may be at most one singular modulus of this discriminant with denominator~$4$. But we already have ${a_z=4}$, and so neither of $a(y^\sigma),a(z^\sigma)$ can equal~$4$ by~\eqref{ealldiff}.

	Thus,  ${a(y^\sigma)=a(z^\sigma)=16}$. Proposition~\ref{plinearel} gives
	$$
	m'+n'\left(1+\frac{1}{4}\right)=\frac{m'}{4}+n'\left(\frac1{16}+\frac1{16}\right),
	$$
	which is  ${m=-3n}$. Arguing as before, 
	we obtain 
	\begin{equation}
	\label{e72}
	\bigl(x^{3}(x^\sigma)^{-3}y^{-1}z^{-1}y^\sigma z^\sigma\bigr)^{24}=1. 
	\end{equation}
	
	We have 
	$$
	\tau(x^\sigma) =\frac{b_4+\sqrt\Delta}{8}, 
	$$
	and we want to specify this $b_4$. Since ${(b_4)^2\equiv \Delta \bmod 16}$, we must have 
	$$
	\Delta \equiv 0,4\bmod 16 \quad\text{and}\quad
	b_4=
	\begin{cases}
	0, &\text{if ${\Delta \equiv 16 \bmod 32}$}, \\
	4, &\text{if ${\Delta \equiv 0 \bmod 32}$},\\
	\pm2, & \text{if ${\Delta \equiv 4 \bmod 16}$}. 
	\end{cases}
	$$
	In particular,
	$$
	b'\in \{b_4+2,b_4-2\}.  
	$$
	Finally, since both $y^\sigma,z^\sigma$ have denominators~$16$ and are $2$-isogenous to~$x^\sigma$, we have 
	$$
	\{\tau(y^\sigma), \tau(z^\sigma)\} =\left\{\frac{b_4+\sqrt\Delta}{16},\frac{b_4'+\sqrt\Delta}{16}\right\} \quad\text{and}\quad
	b_4'\in \{b_4+8,b_4-8\}. 
	$$
	Denote ${t=e^{-\pi|\Delta|^{1/2}/8}}$ and ${\xi=e^{\pi ib_4/8}}$.  Note that ${\xi\in \{1,i, e^{\pm\pi i/4}\}}$, and that 
	$$
	e^{\pi ib'/4}= \pm i\xi^2, \qquad e^{\pi ib_4'/8}=-\xi. 
	$$
	We obtain
	\begin{align*}
	&e^{2\pi i\tau_x}=t^8, \qquad e^{2\pi i \tau_y} =t^{16}, \qquad e^{2\pi i \tau_z} = \eps i \xi^2t^4, \qquad 
	e^{2\pi i \tau(x^\sigma)} = \xi^2t^2,\\
	& \bigl\{e^{2\pi i \tau(y^\sigma)}, e^{2\pi i \tau(z^\sigma)}\bigr\} =\{\xi t,-\xi t\}, 
	\end{align*}
	where ${\eps\in \{1,-1\}}$. 
	Taking the logarithm of~\eqref{e72}, we obtain
	$$
	3 \log(t^8x)-3\log (\xi^2t^2x^\sigma) -\log(t^{16}y)-\log(-\eps i\xi^2t^4z)+\log(-\xi t \cdot \xi t \cdot y^\sigma\cdot z^\sigma) \in \frac{\pi i}{12}\Z. 
	$$
	The $q$-expansion~\eqref{esecond} implies that for some ${k\in \Z}$ we have 
	$$
	\frac{\pi i}{12}k=-162000\xi^2t^2+O_1(10^{10}t^3), 
	$$
	which again leads to a contradiction.

	\subsubsection{Case~\eqref{ecaseonetwoeight}}
	We want to adapt the procedure described in Subsection~\ref{sssfourcases} to this case. We reproduce our data for the reader's convenience:
	\begin{equation}
	\label{esixtwenty}
	\{e_x,e_y\}=\{1,2\},\quad e_z\in\{1,2\}, \quad \Delta\equiv 1\bmod 8,\quad a_x=a_y=1. 
	\end{equation}
	The singular moduli~$x$ and~$y$ are $2$-isogenous by Proposition~\ref{pvarisog}. However, now we have
	${K(x)=K(y)=K(z)=L}$, which means that there does not exist ${\theta \in G}$ with the properties ${x^\theta=x}$ and ${y^\theta =z}$. Hence, a priori we 
	have no control of the degree of isogeny between~$x$ and~$z$. To gain such control we need to determine the denominator~$a_z$. 
	
	\begin{proposition}
		\label{pazletwo}
		Assume~\eqref{esixtwenty}. Then  ${(e_x,e_y)=(1,2)}$ and  
		\begin{equation}
		\label{e1423}
\text{either}\quad		e_z=1, \quad a_z=4 \quad \text{or}\quad 
		e_z=2, \quad a_z= 3. 
		\end{equation}
	\end{proposition}

	The proof consists of several steps. To start with, we eliminate the subcase ${(e_x,e_y)=(2,1)}$. 
	
	\begin{proposition}
		\label{pazleone}
		In case~\eqref{esixtwenty} we must have ${(e_x,e_y)=(1,2)}$.  
	\end{proposition}

	\begin{proof}
		Note that ${a_z>1}$, see~\eqref{eazneone}. 
		We will assume that ${(e_x,e_y)=(2,1)}$ and get a contradiction. 
		
		Since ${\Delta_y=\Delta\equiv 1\bmod 8}$, Proposition~\ref{pevenden} implies that there are~$2$ elements ${\sigma\in G}$ with the property ${a(y^\sigma)=2}$. Since ${L=K(z)}$, at most one of them may satisfies ${a(z^\sigma)=1}$. Hence there exists ${\sigma \in G}$ with the properties 
		${a(y^\sigma)=2}$ and  ${a(z^\sigma) \ge 2}$.

		Since~$x$ and~$y$ are $2$-isogenous, we must have ${a(x^\sigma)\in \{2,8\}}$. But~$2$ is not a denominator for~${\Delta_x=4\Delta}$ by Proposition~\ref{pevenden},
		which implies that ${a(x^\sigma)=8}$. Thus, we found~$\sigma$ such that 
		$$
		a(x^\sigma)=8, \qquad a(y^\sigma)=2, \qquad a(z^\sigma)\ge 2.
		$$
		We now want to derive a contradiction in each of the cases
		\begin{align}
		\label{eone2}
		&\text{one of $a(z),a(z^\sigma)$ is~$2$},\\
		\label{eboth3}
		&\text{both $a(z),a(z^\sigma)\ge3$}. 
		\end{align} 
		
		Assume~\eqref{eone2}. Then ${e_z=1}$, again by the same reason:~$2$ is not a denominator 
		for $4\Delta$. Hence there exists ${\theta\in G}$ such that ${y^\theta=z}$. Since $y,y^\sigma$ are $2$-isogenous, so are ${z=y^\theta}$ and ${z^\sigma=y^{\theta\sigma}=y^{\sigma\theta}}$. It follows that if one of the denominators $a(z),a(z^\sigma)$ is~$2$, then the other must be~$4$. Proposition~\ref{plinearel} now implies that
		$$
		m'+n'+\frac{r'}{a'}= \frac{m'}{8}+\frac{n'}{2}+\frac{r'}{a''} \quad\text{and}\quad \{a',a''\}=\{2,4\}. 
		$$
		Hence
		$$
		\frac78m'= \frac{|n'|}{2}+ r'\left (\frac{1}{a''}-\frac{1}{a'}\right) \le m'\left(\frac12+\frac14\right), 
		$$
		a contradiction. This eliminates~\eqref{eone2}. 
		
		\bigskip

		In case~\eqref{eboth3}
		Proposition~\ref{pineq} used with ${A=9}$ implies that 
		$$
		m' +\frac{|n'|}{2}\le \frac{m'}{8}+ |n'|+\frac{|r'|}{d}+0.01m' \quad\text{and}\quad
		d=\begin{cases}
		\min\{9, a(z^\sigma)\}, & r>0, \\
		\min\{9, a(z)\}, & r<0. 
		\end{cases}
		$$
		Since ${d\ge 3}$, we obtain 
		$$
		\left(\frac78-0.01\right)m' \le \frac{|n'|}{2}+\frac{|r'|}{3} \le m'\left(\frac12+\frac13\right), 
		$$
		a contradiction. This rules~\eqref{eboth3} out as well. The proposition is proved.
	\end{proof}
	
	Next, we show impossibility of ${a_z=2}$.
	
	\begin{proposition}
		\label{paznotwo}
		In case~\eqref{esixtwenty} we must have ${a_z\ge 3}$. 
	\end{proposition}

	\begin{proof}
		We already know that ${a_z\ge2}$ and that ${(e_x,e_y)=(1,2)}$. 
		We also note the statement is immediate for ${e_z=2}$, because~$2$  is not a denominator for $4\Delta$, see Proposition~\ref{pevenden}. Thus, let us assume that 
		$$
		(e_x,e_y,e_z)=(1,2,1), \qquad a_z=2,
		$$
		and show that this is impossible. 
		
		Arguing as in the proof of Proposition~\ref{pazleone} but with the roles of~$x$ and~$y$ interchanged, we find~$\sigma$ satisfying 
		\begin{equation*}
		a(x^\sigma)=2, \qquad a(y^\sigma)=8, \qquad a(z^\sigma)\ge 2. 
		\end{equation*}
		Since $x,z$ are $2$-isogenous, we have ${a(z^\sigma)\in \{1,4\}}$, and we must have ${a(z^\sigma)=4}$ because  ${a(z^\sigma)\ge 2}$. 
		
		Next, let ${\theta\in G}$ be defined by ${z^\theta=x}$. Since $x,z$ are $2$-isogenous, we must have ${a(x^\theta)=2}$, which implies that ${a(y^\theta)=8}$.
		
		Applying Proposition~\ref{plinearel} to the relation
		$$
		(x^\sigma)^m(y^\sigma)^n(z^\sigma)^r= (x^\theta)^m(y^\theta)^n(z^\theta)^r, 
		$$
		we obtain 
		$$
		\frac{m'}{2}+\frac{n'}{8}+\frac{r'}{4}= \frac{m'}{2}+\frac{n'}{8}+\frac{r'}{1}, 
		$$
		which implies ${r=0}$, a contradiction. 
	\end{proof}

	We also need to know that  $|n'|$ is not much smaller than~$m'$.

	\begin{proposition}
		\label{pnprimebelow}
		When ${r>0}$ we have ${|n'|> 0.87 m'}$.
		When ${r<0}$ and ${a_z\ge a}$ we have ${|n'|> \lambda(a) m'}$, where
		$$
		\lambda(a)=0.956-\frac{1}{\min\{30,a\}}. 
		$$
	\end{proposition}
	Here are lower bounds for ${\lambda(a)}$ for some values of~$a$ that will emerge below:
	$$
	\begin{array}{r|rrrrr}
	a&3&5&6&24&30\\
	\hline
	\lambda(a) &>0.62&>0.75&>0.78&>0.91&>0.92
	\end{array}
	$$

	\begin{proof}
		When ${r>0}$ we use Corollary~\ref{cocases} to find~$\sigma$ such that ${a(x^\sigma),a(z^\sigma)\ge 18}$. Now Proposition~\ref{pineq} gives 
		$$
		m' \le |n'| + \frac{m'}{18}+\frac{r'}{18}+0.01 m' \le |n'|+m'\left(\frac2{18}+0.01\right), 
		$$
		which implies ${|n'|> 0.87 m'}$. 
		
		When ${r<0}$ we use Corollary~\ref{cocases} to find~$\sigma$ such that ${a(x^\sigma)\ge 30}$.  When ${a_z\ge a}$, we obtain
		$$
		m'\le |n'|+ \frac{m'}{30}+\frac{|r'|}{\min\{30,a\}} +0.01 m' \le |n'|+m'\left(\frac1{30}+\frac{1}{\min\{30,a\}}+0.01\right),
		$$
		which implies ${|n'|> \lambda(a) m'}$.
	\end{proof}
	
	We are ready to prove Proposition~\ref{pazletwo}.

	\begin{proof}[Proof of Proposition~\ref{pazletwo}]
		The proof is similar to that of Proposition~\ref{pazleone}, but with the roles of~$x$ and~$y$ interchanged. This means that, instead of the  inequality ${m'\ge |n'|}$, we have to use weaker inequalities from Proposition~\ref{pnprimebelow}. This is why we cannot rule out~\eqref{e1423}.

		We already know that ${a_z\ge3}$ and that ${(e_x,e_y)=(1,2)}$. 
		We also note that~$4$ is not a denominator for $4\Delta$, see Proposition~\ref{pevenden}. Hence it suffices to 
		show that each of the cases 
		\begin{align}
		\label{e13}
		&e_z=1, \qquad a_z\ge3, \qquad a_z\ne 4, \\
		\label{e25}
		&e_z=2, \qquad a_z \ge 5
		\end{align}
		leads to a contradiction. As in the proof of Proposition~\ref{paznotwo},  we fix ${\sigma\in G}$ satisfying
		$$
		a(x^\sigma)=2, \qquad a(y^\sigma)=8, \qquad a(z^\sigma)\ge 2. 
		$$
		
		\bigskip
		
		Let us assume~\eqref{e13}. In the same fashion as in the proof of Proposition~\ref{pazleone}, we show that $z,z^\sigma$ are $2$-isogenous. Hence   ${\{a(z),a(z^\sigma)\}=\{a',2a'\}}$, where ${a'\ge 3}$. 
		If ${a'\ge 6}$ then, using 
		Proposition~\ref{pineq}, we obtain 
		$$
		\frac{m'}2 +|n'|\le m'+ \frac{|n'|}8+\frac{|r'|}{6}+0.01m', 
		$$
		which implies that 
		$$
		|n'|\le \frac87m'\left(\frac12+\frac16+0.01\right)<0.78 m', 
		$$ 
		contradicting the lower bound ${|n'|> 0.78m'}$ from  Proposition~\ref{pnprimebelow}.
		
		If ${a'\in \{3,4,5\}}$ then Proposition~\ref{plinearel} gives
		$$
		m'+n'+\frac{r'}{a(z)}=\frac{m'}{2}+\frac{n'}{8}+\frac{r'}{a(z^\sigma)}.
		$$
		This can be rewritten as  
		\begin{equation}
		\label{enprimeid}
		|n'|= \frac87\left(\frac{m'}{2}+r'\left(\frac{1}{a(z)}-\frac{1}{a(z^\sigma)}\right)\right),  
		\end{equation}
		which implies
		\begin{equation}
		\label{e77}
		|n'| \le \frac87m'\left(\frac12+\frac16\right)<0.77m',
		\end{equation}
		When ${r>0}$ this contradicts the lower bound ${|n'|> 0.87m'}$ from Proposition~\ref{pnprimebelow}. When ${r<0}$ and ${a_z=2a'}$  this contradicts  the lower bound ${|n'|> 0.78m'}$ from  Proposition~\ref{pnprimebelow}. Finally, 
		when ${r<0}$ and ${a_z=a'}$, we deduce from~\eqref{enprimeid} the sharper upper bound ${|n'|\le (4/7)m'}$, contradicting 
		the lower bound ${|n'|\ge 0.62m'}$ from Proposition~\ref{pnprimebelow}.   This shows impossibility of~\eqref{e13}.

		\bigskip

		Now let us assume~\eqref{e25}. 
		Proposition~\ref{pineq} implies that 
		$$
		\frac{m'}2 +|n'|\le m'+ \frac{|n'|}8+\frac{|r'|}{d}+0.01m' \quad\text{and}\quad
		d=\begin{cases}
		\min\{9, a(z)\}, & r>0, \\
		\min\{9, a(z^\sigma)\}, & r<0. 
		\end{cases}
		$$
		If ${r>0}$ then ${d\ge 5}$, and we obtain 
		$$
		|n'|\le \frac87m'\left(\frac12+\frac15+0.01\right)<0.82 m', 
		$$
		contradicting the lower bound ${|n'|> 0.87m'}$ from  Proposition~\ref{pnprimebelow}. 
		
		If ${r<0}$ and ${d\ge 8}$ then
		$$
		|n'|\le \frac87m'\left(\frac12+\frac18+0.01\right)<0.73 m', 
		$$ 
		contradicting the lower bound ${|n'|> 0.75m'}$ from  Proposition~\ref{pnprimebelow}.
		
		Thus, 
		$$
		r<0, \qquad 3\le a(z^\sigma)\le 7. 
		$$
		Since ${e_z=2}$, we must have 
		${a(z^\sigma)=p\in \{3,5,7\}}$. Hence~$y^\sigma$ and~$z^\sigma$ are $8p$-isogenous, and so are~$y$ and~$z$. It follows that ${a_z=8p\ge 24}$, and Proposition~\ref{pnprimebelow} implies the lower bound ${|n'|> 0.91m'}$.
		On the other hand, Proposition~\ref{plinearel} implies that 
		$$
		m'+n'+\frac{r'}{8p}= \frac{m'}{2}+\frac{n'}{8}+\frac{r'}{p}, 
		$$
		which yields 
		$$
		|n'|= \frac{4m'}{7}+\frac{|r'|}{p} <0.91m', 
		$$
		a contradiction. This shows impossibility of~\eqref{e25}. 
		The proposition is proved. 
	\end{proof}
	
	Now it is easy to dispose of case~\eqref{ecaseonetwoeight}, alias~\eqref{esixtwenty}. We define~$\sigma_1$ as~$\sigma$ from the proof of Proposition~\ref{pazletwo}; that is, 
	${a(x^{\sigma_1})=2}$ and  ${a(z^{\sigma_1})\ne 1}$. 
	Next, we define~$\sigma_2$ from 
	$$
	z^{\sigma_2}=
	\begin{cases}
	x,& e_z=1,\\
	y,& e_z=2. 
	\end{cases}
	$$ 
	Finally, we set ${\sigma_3=\sigma_2\sigma_1}$.
	
	Using Proposition~\ref{pevenden} and Corollary~\ref{cisog}, we calculate the possible denominators, see Table~\ref{tadensonetwotwoeight}. A verification shows that, in each case, the system of~$3$ linear equation
	\begin{equation}
	\label{esystemthree}
	m'+n' + \frac{r'}{a_z} = \frac{m'}{a(x^{\sigma_i})}+ \frac{n'}{a(y^{\sigma_i})}+ \frac{r'}{a(z^{\sigma_i})} 
	\qquad(i=1,2,3)
	\end{equation}
	has only the trivial solution. 
	
	\begin{table}
		\caption{Denominators for case~\eqref{ecaseonetwoeight} }
		\label{tadensonetwotwoeight}
		{\scriptsize
			$$
			\begin{array}{ccc}
			(e_x,e_y,e_z)=(1,2,1), \quad a_z=4&&
			(e_x,e_y,e_z)=(1,2,2), \quad a_z=3\\
			&&\\
			\begin{array}{cccc}
			i&a(x^{\sigma_i})&a(y^{\sigma_i})&a(z^{\sigma_i})\\
			\hline
			1&2&8&\in \{2,8\}\\
			2&4&16&1\\
			3&\in\{2,8\}&\in \{8,32\}&2
			\end{array}&&
			\begin{array}{cccc}
			i&a(x^{\sigma_i})&a(y^{\sigma_i})&a(z^{\sigma_i})\\
			\hline
			1&2&8&24\\
			2&3&3&1\\
			3&\in \{6,24\}&24&8
			\end{array}\\
			&&\\
			\text{8 systems in total}&&\text{2 systems in total}
			\end{array}
			$$}
	\end{table}

	\subsection{The subdominant case}
	
	\label{sssubdopt}
	
	In this subsection we assume~\eqref{esubdopt}. Thus, we have the following:
	$$
	m>0, \quad n,r<0, \quad m'\ge\max\{ |n'|,|r'|\}, \qquad a_x=1, \quad a_y=a_z=2.
	$$
	To start with, note that 
	\begin{equation}
	\label{eyez}
	e_y=e_z \quad\text{and}\quad \Delta_y=\Delta_z\equiv 1\bmod 8. 
	\end{equation}
	Indeed, among the three numbers $e_x,e_y,e_z$ there are only two distinct integers, see Table~\ref{taexeyez}. If ${e_y\ne e_z}$ then, switching, if necessary,~$x$ and~$y$, we may assume that ${e_x=e_y}$. Hence ${K(x)=K(y)}$, and we have one of the following two possibilities: 
	\begin{align}
	\label{exyzl}
	&K(x)=K(y)=K(z)=L, \\
	&K(x)=K(y)=L, \qquad [L:K(z)]=2. \nonumber
	\end{align}
	In the latter case we must have ${m=n}$ by Theorem~\ref{thprimel}, which is impossible because ${m>0}$ and ${n<0}$.
	Thus, we have~\eqref{exyzl}. Lemma~\ref{ldeq} now implies that ${e_z=2}$ and ${\Delta_x=\Delta_y=\Delta\equiv 1\bmod 8}$. But then ${\Delta_z\equiv 4 \bmod 32}$, and we cannot have ${a_z=2}$ by  Proposition~\ref{pevenden}:\ref{i432none}.
	Thus, ${e_x=e_y}$ is impossible, which proves that ${e_y=e_z}$. Now  Proposition~\ref{pevenden}:\ref{inonemodeight} implies that ${\Delta_y=\Delta_z\equiv 1\bmod 8}$, which completes the proof of~\eqref{eyez}.

	Exploring Table~\ref{taexeyez} and taking note of~\eqref{eyez}, we end up with one of the following cases:
	\begin{align}
	\label{ealll}
	&e_x\in \{1,2\}, && e_y=e_z=1, && \Delta\equiv 1\bmod 8, && L=K(x)=K(y)=K(z);\\
	\label{exylznl}
	&e_x\in\{1,2\},&& e_y=e_z=3, && \Delta\equiv 1 \bmod 24, && [L:K(x)]=2, \quad n=r. 
	\end{align}
Note that we cannot have ${e_y=e_z=4}$ or 	${e_y=e_z=6}$, because in these cases ${\Delta_y=\Delta_z\equiv 0\bmod4}$, contradicting~\eqref{eyez}.

	Each of the cases~\eqref{ealll} and~\eqref{exylznl} can be disposed of using the strategy described in  Subsection~\ref{sssstra}; moreover, the very first step of that strategy can be skipped, because~$a_z$ is already known. 
	
	Case~\eqref{ealll} is analogous to case~\eqref{ecaseonetwoeight}, but is much simpler, because, as indicated above, we already know~$a_z$.  We define $\sigma_1,\sigma_2,\sigma_3$
	by
	$$
	a(y^{\sigma_1})=1, \qquad a(z^{\sigma_2})=1, \qquad a( y^{\sigma_3})=8.  
	$$ 
	Note that there can be several candidates for~$\sigma_3$, we just pick one of them. The possible denominators, determined using Corollary~\ref{cisog} and Proposition~\ref{pevenden}:\ref{i432none},  are given in Table~\ref{tadensalll}. A verification with \textsf{PARI} shows that each of the 12  possible systems 
	$$
	m'+\frac{n'}{2}+ \frac{r'}{2}= \frac{m'}{a(x^{\sigma_i})}+\frac{n'}{a(y^{\sigma_i})}+\frac{r'}{a(z^{\sigma_i})} \qquad(i=1,2,3)
	$$
	has only the trivial solution ${m'=n'=r'=0}$.  
	
	\begin{table}
		\caption{Denominators for case~\eqref{ealll} }
		\label{tadensalll}
		{\scriptsize
			$$
			\begin{array}{ccc}
			(e_x,e_y,e_z)=(1,1,1)&&
			(e_x,e_y,e_z)=(2,1,1) \\
			&&\\
			\begin{array}{cccc}
			i&a(x^{\sigma_i})&a(y^{\sigma_i})&a(z^{\sigma_i})\\
			\hline
			1&2&1&4\\
			2&2&4&1\\
			3&\in\{4,16\}&8&\in\{2,8,32\}
			\end{array}&&
			\begin{array}{cccc}
			i&a(x^{\sigma_i})&a(y^{\sigma_i})&a(z^{\sigma_i})\\
			\hline
			1&8&1&4\\
			2&8&4&1\\
			3&\in \{16,64\}&8&\in\{2,8,32\}
			\end{array}\\
			&&\\
			\text{6 systems in total}&&\text{6 systems in total}
			\end{array}
			$$}
	\end{table}

	\bigskip
	
	In case~\eqref{exylznl} we have ${n=r}$ (and also ${n'=r'}$), and so we need only~$\sigma_1$ and~$\sigma_2$. We do as in Subsection~\ref{sssfourcases}. Since ${\Delta_x\equiv 1 \bmod 3}$, we can find ${\sigma_1,\sigma_2}$ satisfying  
	${a(x^{\sigma_1})=3}$ and ${a( y^{\sigma_2})=9}$.  Defining
	$$
	\ell=
	\begin{cases}
	6, &\text{when ${e_x=1}$}, \\
	12, &\text{when ${e_x=2}$},
	\end{cases}
	$$
	a quick verification shows that singular moduli $x,y$ are $\ell$-isogenous, and so are $x,z$. Using Corollary~\ref{cisog} and   Proposition~\ref{poddpden}:\ref{idivbyp}, we determine the possible denominators: in both cases ${e_x=1}$ and ${e_x=2}$ we find that 
	$$
	a(y^{\sigma_1}), a(z^{\sigma_1})=54, \qquad a(y^{\sigma_2}), a(z^{\sigma_2}) \in\{18,162\}. 
	$$
	It follows that ${m',n'}$ satisfy one of the three linear systems 
	\begin{align*}
	&\begin{cases}
	\displaystyle m'+ n'\left(\frac12+\frac12\right) = \frac{m'}3 + \left(\frac{1}{54}+\frac{1}{54}\right)n', &\\
	&\\
	\displaystyle m'+ n'\left(\frac12+\frac12\right) = \frac{m'}9 +\lambda n', 
	\end{cases}\\
	&\text{where}\quad \lambda \in \left\{\frac{1}{18}+\frac{1}{18}, \frac{1}{18}+ \frac{1}{162}, \frac{1}{162}+\frac{1}{162}\right\}. 
	\end{align*}
	A verification shows that each of these systems has only the trivial solution ${m'=n'=0}$. This completes the proof of Theorem~\ref{thpizthree} for equal fundamental discriminants.

	\subsection{Distinct fundamental discriminants}
	\label{ssdis}
	We are left with the case when the fundamental discriminants  $D_x,D_y,D_z$ are not all equal. We may assume that
	${|\Delta_z|\ge |\Delta_y| \ge |\Delta_x|}$. In particular, ${|\Delta_z|\ge 10^{10}}$ and ${\Delta_z\ne \Delta_x}$. 
	Theorem~\ref{thprimel} and Lemma~\ref{lnorou} imply that 
	\begin{equation*}
	\Q(y)=\Q(y^n)= \Q(x^mz^r) =\Q(x,z). 
	\end{equation*}
	In particular, 
	\begin{equation}
	\label{exziny}
	\Q(x), \Q(z)\subset \Q(y). 
	\end{equation}
Furthermore, Theorem~\ref{thprimel} and Lemma~\ref{lnorou} imply that ${\Q(x)=\Q(x^m)=\Q(y^nz^r)}$ is a subfield of $\Q(y,z)$ of degree at most~$2$. Since ${z\in \Q(y)}$, this implies that 
$$
[\Q(y):\Q(x)]\le 2.
$$  	
Unfortunately, we cannot claim similarly that ${[\Q(y):\Q(z)]\le 2}$, because we do not know whether the singular moduli $x,y$ satisfy the hypothesis of Theorem~\ref{thprimel}. 

The rest of the proof splits into two cases: ${D_x\ne D_y}$ and ${D_y\ne D_z}$.

	\subsubsection{The case ${D_x\ne D_y}$}
Since 	${[\Q(y):\Q(x)]\le 2}$,  Corollary~\ref{cxy}  implies that  ${[\Q(y):\Q]\le 32}$. It follows that  ${h(\Delta_z)=[\Q(z):\Q]\le 32}$ as well. This contradicts Proposition~\ref{pwatki} because 
${|\Delta_z|\ge 10^{10}}$. 

\subsubsection{The case ${D_y\ne D_z}$}
Since ${z\in \Q(y)}$,  the field $\Q(z)$ must be $2$-elementary by Proposition~\ref{pintersect}. Hence the proof in this case will be complete if we show one of the following:
\begin{align}
\label{erhotwoz}
\rho_2(\Delta_z) &\le 6;\\
\label{eqyqz}
[\Q(y):\Q(z)]&\le 2. 
\end{align}
Recall that $\rho_2(\cdot)$ is the $2$-rank, see Subsection~\ref{ssstwor}.

Indeed, assume that~\eqref{erhotwoz} holds.  Since $\Q(z)$ is $2$-elementary, we have ${h(\Delta_z) = 2^{\rho_2(\Delta_z)}\le 64}$,  contradicting Proposition~\ref{pwatki}.

Similarly, if~\eqref{eqyqz} holds then ${h(\Delta_z)\le 16}$ by Corollary~\ref{cxy},  again contradicting Proposition~\ref{pwatki}.

We will show that, depending on the size of~$\Delta_y$,  one of~\eqref{erhotwoz} or~\eqref{eqyqz} holds. 

\bigskip
	
If
 ${|\Delta_y|\ge10^8}$, then Theorem~\ref{thprimel} applies to the singular moduli $x,y$. It follows that 
${\Q(z)=\Q(z^r)=\Q(x^my^n)}$ is a subfield of ${\Q(x,y)}$ of degree at most~$2$. Since ${x\in \Q(y)}$, we obtain 
${[\Q(y):\Q(z)]\le 2}$, which is~\eqref{eqyqz}.  
	
\bigskip

%\subsubsection{The subcase ${10^6\le |\Delta_y|\le 10^8}$, ${D_y\ne D_z}$}
Next, let us assume that ${10^6\le |\Delta_y|\le 10^8}$. 
If ${\Delta_y\not\equiv 4\bmod 32}$ then Theorem~\ref{thprimel} again applies to the singular moduli $x,y$, and we may argue as above, obtaining~\eqref{eqyqz}.  

If ${\Delta_y\equiv 4\bmod 32}$, then ${\rho_2(\Delta_y)= \omega(\Delta_y)-2}$, see Proposition~\ref{pgauss}. Since 
$$
10^8<4\cdot 3\cdot5\cdot7\cdot11\cdot13\cdot17\cdot 19\cdot 23=446185740,
$$
we have ${\omega(\Delta_y) \le 8}$. Hence ${\rho_2(\Delta_y)\le6}$. 
Now let ${K=\Q(\sqrt{\Delta_y})}$ be the CM-field of~$y$. Since $\Q(z)$ is $2$-elementary,~$K$ is not contained in $\Q(z)$ by Proposition~\ref{preal}.  Since both~$K$ and~$\Q(z)$ are Galois extensions of~$\Q$,  the group $\Gal(\Q(z)/\Q)$  is isomorphic to $\Gal(K(z)/K)$, which is a quotient group of $\Gal(K(y)/K)$. In particular,
$$
\rho_2(\Gal(\Q(z)/\Q))\le \rho_2(\Gal(K(y)/K))\le 6,
$$
which is~\eqref{erhotwoz}.

\bigskip

%\subsubsection{The subcase , ${D_y\ne D_z}$}	
Finally, let us assume that ${|\Delta_y|\le 10^6}$. 
Since
$$
10^6< 4\cdot 3\cdot 5\cdot 7\cdot 11\cdot 13 \cdot 17= 1021020,
$$ 
we have
%Hence every~$\Delta_y$ with ${|\Delta_y|\le 10^6}$ satisfies 
${\rho_2(\Delta_y)\le \omega(\Delta) \le 6}$, again by Proposition~\ref{pgauss}.  As we have seen above, this implies~\eqref{erhotwoz}. 
Theorem~\ref{thpizthree} is proved.

	{\footnotesize

		\bibliographystyle{amsplain}
		\bibliography{pizthree}

\noindent
\textbf{Yuri Bilu \& Emanuele Tron}\\
Institut de Mathématiques de Bordeaux\\
Université de Bordeaux \& CNRS\\
33405 Talence France\\
\url{yuri@math.u-bordeaux.fr}\\
\url{emanuele.tron@math.u-bordeaux.fr}

\bigskip

\noindent
\textbf{Sanoli Gun}\\
The Institute of Mathematical Sciences\\
A CI of Homi Bhabha National Institute \\
CIT Campus, Taramani, Chennai 600 113 \\
India\\
\url{sanoli@imsc.res.in}

}
\end{document}